\documentstyle[12pt]{article}
\newtheorem{theorem}{{\sc Theorem}}
\newcommand{\bt}{\begin{theorem}}
\newcommand{\et}{\end{theorem}}
\newcommand{\newsection}[1]{\setcounter{equation}{0} \setcounter{theorem}{0}
\section{#1}}

\newcommand{\NI}{\noindent}

\newcommand{\bea}{\begin{eqnarray}}
\newcommand{\eea}{\end{eqnarray}}

\def \b #1 {\bf #1}
\newcommand{\IR}{I\!\!R}
\newcommand{\IE}{I\!\!E}
\newcommand{\IC}{I\!\!C}

\newcommand{\IT}{I\!\!T}
\newcommand{\IN}{I\!\!N}
\newcommand{\IZ}{Z\!\!\!Z}
\newcommand{\clk}{{\cal K}}

\newcommand{\cla}{{\cal A}}

\newcommand{\cli}{{\cal I}}

\newcommand{\clf}{{\cal F}}
\newcommand{\clg}{{\cal G}}
\newcommand{\clh}{{\cal H}}
\newcommand{\clp}{{\cal P}}
\newcommand{\clo}{{\cal O}}
\newcommand{\clb}{{\cal B}}

\newcommand{\cle}{{\cal E}}

\newcommand{\cln}{{\cal N}}

\newcommand{\clm}{{\cal M}}

\newcommand{\raro}{\rightarrow}

\newcommand{\vsp}{\vskip 1em}

\def \qed {\hfill \vrule height6pt width 6pt depth 0pt}
\newcommand{\be}{\begin{equation}}
\newcommand{\ee}{\end{equation}}
\newcommand{\ben}{\begin{eqnarray*}}
\newcommand{\een}{\end{eqnarray*}}
\begin{document}
\thispagestyle {empty}
\sloppy

\centerline{\large \bf Markov shift in non-commutative probability -II }

\bigskip
\centerline{\bf Anilesh Mohari }
\smallskip
\centerline{\bf S.N.Bose Center for Basic Sciences, }
\centerline{\bf JD Block, Sector-3, Calcutta-98 }
\centerline{\bf E-mail:anilesh@boson.bose.res.in}
\smallskip
\centerline{\bf Abstract}
\bigskip
We study asymptotic behavior of a Markov semigroup on a von-Neumann algebra by 
exploring a maximal von-Neumann subalgebra where the Markov semigroup is an 
automorphism. This enables us to prove that strong mixing is equivalent to ergodic property 
for continuous time Markov semigroup on a type-I von-Neumann algebra with center 
completely atomic. For discrete time dynamics we prove that an aperiodic ergodic Markov 
semigroup on a type-I von-Neumann algebra with center completely atomic is strong mixing. There exists 
a tower of isomorphic von-Neumann algebras generated by the weak Markov process
and a unique up to isomorphism minimal dilated quantum dynamics of endomorphisms associated with the Markov
semigroup. The dilated endomorphism is pure in the sense of Powers if and only if the adjoint Markov semigroup 
satisfies Kolmogorov property. As an application of our general results we find a necessary and sufficient condition 
for a translation invariant state on a quantum spin chain to be pure. We also find a tower of type-II$_1$ factors 
canonically associated with the canonical conditional expectation 
on a sub-factor of a type-II$_1$ factor. This tower of factors unlike Jones's tower do not preserve index. This
gives a sequence of Jones's numbers as an invariance for the inclusion of a finite sub-factor of a type-II$_1$ 
factor.

\newpage
\newsection{ Introduction:}

\vsp
Let $\tau=(\tau_t,\;t \ge 0)$ be a semigroup of identity preserving completely
positive normal maps [Da,BR] on a von-Neumann algebra $\cla_0$ acting on a
separable Hilbert space $\clh_0$, where either the parameter $t \in \!R_+$, the
set of positive real numbers or $\!Z^+$, the set of positive integers. In case $t \in \!R_+$, i.e.
continuous, we assume that for each $x \in \cla_0$ the map $t \raro \tau_t(x)$ is continuous in 
the weak$^*$ topology. Thus variable $t \in \IT_+$ where $\IT$ is either $\IR$ or $\IN$. We assume further 
that $(\tau_t)$ admits a normal invariant state $\phi_0$, i.e. $\phi_0 \tau_t = \phi_0 \forall t \ge 0$. 
We continue our investigation [Mo1] on asymptotic behavior of the quantum dynamical semigroup $(\cla_0,\tau_t,\phi_0)$ 
and associated minimal dilated processes.

\vsp
In section 2 we investigate asymptotic behavior of the quantum dynamical semigroup $(\tau_t\;t \ge 0)$ on $\cla_0$. We say $(\tau_t)$ is {\it ergodic } if $\{ x: \tau_t(x)=x, \; t \ge 0 \}= \{ zI, z \in \IC \}$ and a 
normal state $\phi_0$ is {\it invariant} if $\phi_0(\tau_t(x))=\phi_0(x)$ for all $x \in \cla_0, t \ge 0$. A normal state $\phi_0$ is
an {\it equilibrium or strongly mixing } state if $\phi \tau_t(x) \raro  \phi_0(x)$ as $t \raro \infty$ for all $x \in \cla_0$
and normal state $\phi$ on $\cla_0$. Let $p$ be the support projection of the state $\phi_0$ in $\cla_0$. $p$ is the minimal 
element in $\cla_0$ so that $\phi_0(x)=\phi_0(pxp)$. Thus $p$ is a {\it sub-harmonic} projection (i.e. $\tau_t(p) \ge p$ ) for 
$(\tau_t)$ and $\tau_t(x)p=\tau_t(xp)p$ for all $x \in \cla_0$. We define the {\it reduced Markov semigroup} 
$(\cla_0^p,\tau^p_t,\phi^p_0)$ by $\tau^p_t(x)=p\tau_t(pxp)p$ for all $x \in \cla^p_0$, where $\cla^p_0=p\cla_0p$. 
Let $y$ be the strong limit of $\tau_t(p)$ as $t \raro \infty$. So $\tau_t(y)=y$ for all $t \ge 0$. Thus $(\cla_0,\tau_t,\phi_0)$ 
is ergodic if and only if $\tau_t(p) \uparrow I$ and $(\cla_0^p,\tau^p_t,\phi^p_0)$ is ergodic ( See Theorem 3.6 in [Mo1]). Similar 
result is also for strong mixing ( see Theorem 3.12 in [Mo1] ). All these results suggest while studying asymptotic behavior, there 
is no charm lost in assuming that the invariant normal state is faithful. In this section we aim to refine various sufficient conditions 
proved in [Mo1] for strong mixing.         

\vsp
In case $\phi_0$ is faithful, normal and invariant for $(\tau_t)$, we recall [Mo1] that
$\clg=\{x \in \cla_0: \tilde{\tau}_t \tau_t(x) = x,\; t \ge 0 \}$ is von-Neumann sub-algebra of $\clf= \{ x  \in \cla_0:
\tau_t(x^*)\tau_t(x) = \tau_t(x^*x),\tau_t(x)\tau_t(x^*) = \tau_t(xx^*)\; \forall t \ge 0 \} $ and
the equality $\clg=\IC$ is a sufficient condition for $\phi_0$ to be strong mixing for $(\tau_t)$. Since the
backward process [AM] is related with the forward process via an anti-unitary operator we note that $\phi_0$ is
strongly mixing for $(\tau_t)$ if and only if same hold for $(\tilde{\tau}_t)$. We can also
check this fact by exploring faithfulness of $\phi_0$ and the adjoint relation [OP]. Thus $\IC \subseteq \tilde{\clg} \subseteq
\tilde{\clf}$ and equality $\IC=\tilde{\clg}$ is also a sufficient condition for strong mixing where $\tilde{\clf}$ and
$\tilde{\clg}$ are von-Neumann algebras associated with $(\tilde{\tau}_t)$. Thus we find two competing criteria for strong mixing.
However it is not clear whether $\clf = \tilde{\clf}$ or $\clg=\tilde{\clg}$. Since given a dynamics it is difficult to describe
$(\tilde{\tau}_t)$ explicitly this criterion $\clg=\IC$ is rather non-transparent. We prove that
$\clg=\{ x \in \clf: \tau_t\sigma_s(x) = \sigma_s \tau_t(x),\; \forall t \ge 0.\; s \in \!R \}$ where
$\sigma=(\sigma_s:\;s \in \!R)$ is the Tomita's modular auto-morphism group [BR,OP] associated with $\phi_0$.
So $\clg$ is the maximal von-Neumann sub-algebra of $\cla_0$, where $(\tau_t)$ is an $*$-endomorphism [Ar], invariant by the
modular auto-morphism group $(\sigma_s)$. 
Moreover $\sigma_s(\clg)=\clg$ for all $s \in \!R$ and $\tilde{\tau}_t(\clg)=\clg$ for all $t \ge 0$.
Thus by a theorem of Takesaki [OP], there exists a norm one projection $\IE_{\clg}$
from $\cla_0$ onto $\clg$ which preserves $\phi_0$ i.e. $\phi_0 \IE = \phi_0$. Exploring
the fact that $\tilde{\tau}_t(\clg)=\clg$, we also conclude that the conditional expectation $\IE_{\clg}$ commutes with
$(\tau_t)$. This enables us to prove that $(\cla_0,\tau_t,\phi_0)$ is ergodic (strongly mixing) if and only if $(\clg,\tau_t,\phi_0)$
is ergodic (strongly mixing). Though $\tau_t(\clg) \subseteq \clg$ for all $t \ge 0$, equality may not hold in general.
However we have $$\bigcap_{t \ge 0} \tau_t(\clg)=\bigcap_{t \ge 0}\tilde{\tau}_t(\tilde{\clg})$$
where $\tilde{\clg}=\{ x \in \cla_0: \tau_t(\tilde{\tau}_t(x))=x,\; t \ge 0 \}$.
$\clg = \tilde{\clg}$ holds if and only if $\tau_t(\clg)=\clg,\;\tilde{\tau}_t(\tilde{\clg})=\tilde{\clg}$ for all $t \ge 0$. 
Thus $\clg_0= \bigcap_{t \ge 0} \tau_t(\clg)$ is the maximal von-Neumann sub-algebra invariant by the modular
automorphism so that $(\clg_0,\tau_t,\phi_0)$ is an $*-$automorphisms with $(\clg_0,\tilde{\tau}_t,\phi_0)$
as it's inverse dynamics. Once more there exists a conditional expectation $\IE_{\clg_0}:\cla_0 \raro \cla_0$
onto $\clg_0$ commuting with $(\tau_t)$. This ensures that $(\cla_0,\tau_t,\phi_0)$ is ergodic (strongly mixing)
if and only if $(\clg_0,\tau_t,\phi_0)$ is ergodic (strongly mixing). It is clear now that
$\clg_0=\tilde{\clg}_0$, thus $\clg_0=\IC$, a criterion for strong mixing, is symmetric or time-reversible.
Exploring the criterion $\clg_0=\IC$ we also prove that for a type-I factor $\cla_0$ with center completely atomic,
strong mixing is equivalent to ergodicity when the time variable is continuous i.e. $\!R_+$ (Theorem 2.4). This result in particular
extends a result proved by Arveson [Ar] for type-I finite factor. In general, for discreet time dynamics $(\cla_0,\tau,\phi_0)$, 
ergodicity does not imply strong mixing property (not a surprise fact since we have many classical cases). We prove that 
$\tau$ on a type-I von-Neumann algebra $\cla_0$ with completely atomic center is strong mixing 
if and only if it is ergodic and the point spectrum of $\tau$ in the unit circle i.e. $\{w \in S^1: \tau(x)=wx 
\;\;\mbox{for some non zero}\;\; x in \cla_0 \}$ is trivial.

\vsp
In section 3 we consider the unique up to isomorphism minimal forward weak Markov [AM,Mo1] stationary process $\{ 
j_t(x),\; t \in \IT,\; x \in \cla_0 \}$ associated with $(\cla_0,\tau_t,\phi_0)$. We set a family of isomorphic von-Neumann 
algebras $\{ \cla_{[t}: t \in \IT \}$ generated by the forward process so that $\cla_{[t} \subseteq \cla_{[s}$ whenever 
$s \le t$. In this framework we construct a unique modulo unitary equivalence minimal 
dilation $(\cla_{[0},\alpha_t,\; t \ge 0,\phi)$, where $\alpha=(\alpha_t:t \ge 0)$ is a semigroup of $*-$endomorphism on 
a von-Neumann algebra $\cla_{[0}$ acting on a Hilbert space $\clh_{[0}$ with a normal invariant state $\phi$ and a 
projection $P$ in $\cla_{[0}$ so that

\NI (a) $P\cla_{[0}P=\pi(\cla_0)''$; 

\NI (b) $\Omega \in \clh_{[0}$ is a unit vector so that $\phi(X)=<\Omega, X \Omega>$;

\NI (b) $P\alpha_t(X)P=\pi(\tau_t(PXP))$ for $t \ge 0,\;X \in \cla_{[0}$; 

\NI (c) $ \{ \alpha_{t_n}(PX_nP) .....\alpha_{t_3}(PX_3P)\alpha_{t_2}(PX_2P)\alpha_{t_1}(PX_1P)\Omega:\; 0 \le t_1 \le t_2 ..\le t_n
,\; n \ge 1 \}, X_i \in \cla_{[0} \}$ is total in $\clh_{[0}$,

\NI where $\pi$ is the GNS representation of $\cla_0$ associated with the state $\phi_0$. We end this section with a criterion 
for the inductive limit state associated with a $C^*$ algebra valued quantum dynamical semigroup 
$(\clb_0,\lambda_t:t \ge 0,\psi)$ of endomorphisms to be pure. To that end we explore the minimal weak Markov process associated
with the reduced Markov semigroup on the corner algebra of the support projection and prove that the inductive limit state
is pure if and only if the Markov semigroup satisfies Kolmogorov's property introduced in [Mo1].

\vsp
In section 4 we deal with only faithful $\psi_0$ and prove that $\cla_{[t}$ is a factor if and only if $\cla_0$ is a factor.
Moreover $\cla_{[t}$ is a type-I (type-II, type-III) factor if and only if $\cla_0$ is also type-I (type-II, type- III) 
respectively. In particular we construct a class of complete boolean algebra of factors [AW]. In section 5 we find a product 
system [Ar] when $\cla_0$ is a type-I factor and construct a class of complete boolean algebra of type-I factors [AW] 
appearing canonically with a continuous tensor product of Hilbert spaces.

\vsp
Section 6 includes an application of results proved in section 3. We consider a translation invariant extremal 
state $\omega'$ on UHF$_d$ algebra $\otimes_{\IZ} M_d$ and choose an element $\psi \in K_{\omega}$, where $K_{\omega}= 
\{ \psi: \psi \mbox{ is a state on } \clo_d \mbox{ such that } \psi \lambda = \psi \mbox{ and } \psi_{|\mbox{UHF}_d}= \omega \}$, 
$\lambda$ is the canonical endomorphism on $\clo_d$ and $\omega$ is the restriction of $\omega'$. Let $(\clh_{\pi},\pi,\Omega)$ 
be the GNS representation of $(\clo_d,\psi)$. Then $(\clh,S_i,P,V_i,\Omega)$ is a Popescu system [BJKW], where 
$P$ is the support projection of the state $\psi_{\Omega}(X)=<\Omega,X\Omega>$ on the von-Neumann algebra $\pi(\clo_d)''$,
$S_i=\pi(s_i)$ and $V_i=PS_iP$. Let $Q$ be the support projection of the state $\psi_{\Omega}$ on the von-Neumann
algebra $\{S_IS^*_J: |I|=|J| < \infty  \}''$ and $\cla_0=QS_IS^*_JQ: \; |I|=|J| < \infty \}$. We also 
set $l_k=QS_kQ$ for all $ 1 \le k \le d$ and define Markov semigroup $\tau$ on $\cla_0$ by
$\tau(x)=\sum_i l_ixl_i^*$. The normal state $\psi_0$, defined by $\psi_0(x) =\psi(QxQ)$ for all $x \in \cla_0$, is
faithful normal and invariant for $\tau$. We explore Kolmogorov's property of the minimal weak Markov process associated 
with $(\cla_0,\tau_n,\psi_0)$ is a necessary and sufficient condition for $\omega'$ to be pure. The result here is 
more general then what initiated and developed in [FNW1,FNW2, BJKW] for translational invariant state on quantum spin chain. 
The theory is further developed when the state is in detailed balance [Mo3] and applied to study behaviour of the ground
state for well known examples.

\vsp
In section 7 we investigate the tower $\cla_{[t}$ of factors when $\cla_0$ is type-II$_1$. In such a case each $\cla_{[t}$ is either 
identical and isomorphic to $\cla_0$ or is a type-II$_\infty$ factor. Moreover $j_0(I)$ is a finite projection in $\cla_{[-t}$ for each 
$t \ge 0$. Thus we find a canonical tower of type-II$_1$ factors $\clm_s \subseteq \clm_t$ for $s \le t$, where $\clm_t=j_0(I)\cla_{[-t}
j_0(I),\; t \ge 0$, acting on the Hilbert subspace $j_0(I)$. One natural question that appears interesting: How Jones's tower of 
type-II$_1$ factors is related with the tower $\{ \clm_t : t \ge 0 \}$? Can we recover Jones's tower by choosing an appropriate 
dynamics $(\cla_0,\tau,\phi_0)$ in discreet time variable? To that end let $\clb_0$ be a proper finite sub-factor of $\cla_0$ and 
$\phi_0$ be the unique normalize trace. We consider the representation of $\cla_0$ by left multiplication on $L^2(\cla_0,\phi_0)$ and 
a conditional expectation $\tau$ on $\clb_0$ defined by $\tau(x)=E_0xE_0$, where $E_0$ is the projection in $L^2(\cla_0,\phi_0)$ 
generated by vectors in $\clb_0$. We prove that $\cla_1 =\{ \cla_0, E_0 \}''$ is isomorphic to a proper von-Neumann sub-algebra of 
$\clm_1$ associated with $(\cla_0,\tau, \phi_0)$. Thus $[\clm_1: \clm_0] > [\cla_0 : \clb_0]$. So the canonical tower $\clm_k \subseteq 
\clm_{k+1} $ of type-II$_1$ factors appears here is different from that of Jones's [Jo]. The sequence $\{ [\clm_k : \clm_{k-1}]:
\;\; k \ge 0 \}$ of Jones index is an invariance for the inclusion of the sub-factors. A detailed study, needs to be done to explore 
this new invariance, which seems to be an interesting problem!

\newsection{ Time-reverse Markov semigroup and asymptotic properties: }

\vsp
Following [OP,AM], we consider the unique Markov map $\tilde{\tau}$ on $\cla_0$ 
which satisfies the following adjoint relation
\be
\phi_0(\sigma_{1/2}(x)\tau(y))=\phi_0(\tilde{\tau}(x)\sigma_{-1/2}(y))
\ee
for all $x,y \in \cla_0$ analytic elements for the Tomita's modular 
automorphism $(\sigma_t:\; t \in \IR)$ associated with a faithful normal 
invariant state for a Markov map $\tau$ on $\cla_0$. For more details 
we refer to the monograph [OP]. We also quote now [OP, Proposition 8.4 ] the 
following proposition without a proof. 

\vsp
\NI {\bf PROPOSITION 2.1: } Let $\tau$ be an unital completely positive normal
maps on a von-Neumann algebra $\cla_0$ and $\phi_0$ be a faithful normal 
invariant state for $\tau$. Then the following conditions are equivalent 
for $x \in \cla_0$:

\NI (a) $\tau(x^*x)=\tau(x^*)\tau(x)$ and $\sigma_s(\tau(x))= \tau(\sigma_s(x)),\; \forall \; s \in \!R;$

\NI (b) $\tilde{\tau} \tau (x)=x.$  

\NI Moreover $\tau$ restricted to the sub-algebra $\{x: \tilde{\tau}\tau(x) 
=x \}$ is an isomorphism onto the sub-algebra $\{x \in \cla_0: \tau\tilde{\tau}
(x) =x \}$ where $(\sigma_s)$ be the modular automorphism on $\cla_0$ 
associated with $\phi_0$. 

\vsp
\NI {\bf PROPOSITION 2.2: } Let $(\cla_0,\tau_t,\phi_0)$ be a quantum dynamical system and $\phi_0$ be 
faithful invariant normal state for $(\tau_t)$. Then the following hold:

\NI (a) $\clg = \{x \in \cla_0: \tau_t(x^*x)=\tau_t(x^*) \tau_t(x),\; \tau_t(xx^*) = \tau_t(x)\tau_t(x^*),\; \sigma_s(\tau_t(x)) 
        = \tau_t (\sigma_s(x)),\; \forall \; s \in \!R,\; t \ge 0 \}$ and $\clg$ is $\sigma= (\sigma_s:\;s \in \!R)$ invariant 
        and commuting with $\tau=(\tau_t:t \ge 0)$ on $\clg$. Moreover for all $t \ge 0,\; \tilde{\tau}_t(\clg)=\clg$ and the 
        conditional expectation $\!E_{\clg}: \cla_0 \raro \cla_0$ onto $\clg_0$ commutes with $(\tau_t)$.

\NI (b) There exists a unique maximal von-Neumann algebra $\clg_0 \subseteq \clg \bigcap \tilde{\clg}$ so that $\sigma_t(\clg_0)=\clg_0$ 
for all $t \in \!R$ and $(\clg_0,\tau_t,\phi_0)$ is an automorphism where for any $t \ge 0$, $\tilde{\tau}_t\tau_t=\tau_t\tilde{\tau}_t=1$ 
on $\clg_0$. Moreover the conditional expectation $\!E_{\clg_0}:\cla_0 \raro \cla_0$ onto $\clg_0$ commutes with $(\tau_t)$ and $(\tilde{\tau}_t)$.

\vsp
\NI {\bf PROOF:} The first part of (a) is a trivial consequence of Proposition 2.1 
once we note that $\clg$ is closed under the action $x \raro x^*$. For the 
second part we recall [Mo1] that $\phi_0(x^*JxJ)- 
\phi_0(\tau_t(x^*)J\tau_t(x)J)$ is monotonically increasing with $t$ and 
thus for each $t \ge 0$ if $\tilde{\tau}_t\tau_t(x)=x$ then 
$\tilde{\tau}_s\tau_s(x)=x$ for all $0 \le s \le t$. So the 
sequence $\clg_t=\{ x \in \cla_0:\tilde{\tau}_t\tau_t(x)=x \}$ of von-Neumann 
sub-algebras decreases to $\clg$ as $t$ increases to $\infty$ i.e.
$\clg = \bigcap_{t \ge 0} \clg_t$. Similarly we also have $\tilde{\clg}=
\bigcap_{ t \ge 0} \tilde{\clg}_t$.   

\vsp
Since $\tilde{\clg}_t$ monotonically decreases to $\tilde{\clg}$ as $t $ increases to infinity 
for any $s \ge 0$ we have $\tau_s(\tilde{\clg}) = \bigcap_{t \ge 0} \tau_s(\tilde{\clg}_t)$.

Now we verify that $\bigcap_{s \ge r } \tau_s(\tilde{\clg}) = \bigcap_{ s \ge r } \bigcap_{t \ge 0} \tau_{s+t}(\clg_t)=
\bigcap_{t \ge 0} \bigcap_{s \ge r} \tau_{s+t}(\clg_t) = \bigcap_{t \ge r} \bigcap_{0 \le s \le t}\tau_t(\clg_s)$, 
where we have used $\tau_t(\clg_t)=\tilde{\clg}_t$. Since $\clg_t$ 
are monotonically decreasing with $t$ we also note that 
$\bigcap_{0 \le s \le t} \tau_t(\clg_s) = \tau_t(\clg_t)$. Hence for any $r \ge 0$
\be
\bigcap_{s \ge r } \tau_s(\tilde{\clg}) = \tilde{\clg} 
\ee
From (2.2) with $r=0$ we get $\tilde{\clg} \subseteq \tau_t(\tilde{\clg})$ for all $t \ge 0$. For any $t \ge 0$ 
we also have $\tau_t(\tilde{\clg}) \subseteq \bigcap_{s \ge t } \tau_s(\tilde{\clg}) = \tilde{\clg}$. Hence we
conclude $\tau_t(\tilde{\clg}) = \tilde{\clg}$ for any $t \ge 0$. By symmetry $\tilde{\tau}_t(\clg)=\clg$ for any $t \ge 0$.  

\vsp
Since $\clg$ is invariant under the modular automorphism $(\sigma_s)$ by a theorem of
Takesaki [AC] there exists a norm one projection $\!E_{\clg}:\cla \raro \cla$ with range equal
to $\clg$. We claim that $\!E_{\clg}$ commutes with $(\tau_t)$. To that end we verify for any $x \in \cla_0$ and
$y \in \clg$ the following equalities:

$$ <J_{\clg}yJ_{\clg}\omega_0, \!E_{\clg}(\tau_t(x)) \omega_0> = <J_0yJ_0\omega_0, \tau_t(x) \omega_0>$$
$$=<J_0\tilde{\tau}_t(y)J_0 \omega_0,x \omega_0> = <J_{\clg}\tilde{\tau}_t(y)J_{\clg} \omega_0, \!E_{\clg}(x)
\omega_0> $$
$$ = <J_{\clg}yJ_{\clg}\omega_0,\tau_t(\!E_{\clg}(x))\omega_0)> $$
where we used the fact that $\tilde{\tau}(\clg)=\clg$ for the third equality and range of $\IE_{\clg}$ is indeed $\clg$ is
used for the last equality. This completes the proof of (a).

\vsp
Now for any $s \ge 0$ we note that 
$\tilde{\tau}_s(\tilde{\clg}) = \bigcap_{t \ge s} \tilde{\tau}_s(\tilde{\clg}_t) = \bigcap_{t \ge s} \tilde{\tau}_s 
\tau_t(\clg_t)= \bigcap_{t \ge s}\tilde{\tau}_s \tau_s(\tau_{t-s}(\clg_t)))=\bigcap_{t \ge s} \tau_{t-s}(\clg_t)=
\bigcap_{t \ge 0}\tau_t(\clg_{s+t})$ where we have used the fact that $\tau_{t-s}(\clg_t) \subseteq \clg_s$. 
Thus we have $\bigcap_{s \ge 0}\tau_s(\clg) \subseteq  \bigcap_{s \ge 0}\tilde{\tau}_s(\tilde{\clg})$. By the dual 
symmetry, we conclude the reverse inclusion and hence  
\be
\bigcap_{s \ge 0}\tau_s(\clg) =  \bigcap_{s \ge 0}\tilde{\tau}_s(\tilde{\clg})
\ee

\vsp
We set $\clg_0=\bigcap_{s \ge 0}\tau_s(\clg)$. Thus $\clg_0 \subseteq \clg$ and also $\clg_0 \subseteq \tilde{\clg}$ by (2.3) and 
for each $t \ge 0$ we have $\tau_t\tilde{\tau}_t=\tilde{\tau}_t \tau_t=1$ on $\clg_0$. Since $\tau_s(\clg)$ is monotonically decreasing, we also
note that $\tau_t(\clg_0)= \bigcap_{s \ge 0}\tau_{s+t}(\clg)=\clg_0$. Similarly $\tilde{\tau}_t(\clg_0)=\clg_0$ by (2.3). 
That $\clg_0$ is invariant by the modular group $\sigma$ follows since $\clg$ is invariant by $\sigma=(\sigma_t)$ which is commuting 
with $\tau=(\tau_t)$ on $\clg$. Same is also true for $(\tilde{\tau}_t)$ by (2.3). By Takesaki's theorem [AC] once more we guarantee
that there exists a conditional expectation $\!E_{\clg_0}: \cla_0 \raro \cla_0$ with range equal to $\clg_0$. Since $\tilde{\tau_t}(\clg_0)=
\clg_0$, once more by repeating the above argument we conclude that $E_{\clg_0} \tau_t= \tau_t E_{\clg_0}$ on $\cla_0$. Since we also have 
$\tau_t(\clg_0)=\clg_0$, by symmetry of the argument $\!E_{\clg_0}$ is also commuting with $\tilde{\tau}=(\tilde{\tau}_t)$ \qed

\vsp
We have the following reduction theorem.  

\vsp
\NI {\bf THEOREM 2.3: } Let $(\cla_0,\tau_t,\phi_0)$ be as in Proposition 2.2. 
Then the following statements are equivalent:

\NI (a) $(\cla_0,\tau_t,\phi_0)$ is strong mixing ( ergodic );

\NI (b) $(\clg, \tau_t,\phi_0)$ is strong mixing ( ergodic );

\NI (c) $(\clg_0,\tau_t,\phi_0)$ is strong mixing ( ergodic ). 

\vsp
\NI {\bf PROOF: } That (a) implies (b) is obvious. By Proposition 2.2. we have $\!E_{\clg}\tau_t(x) = \tau_t \!E_{\clg}(x)$ for any 
$x \in \cla_0$ and $t \ge 0$. Fix any $x \in \cla_0$. Let $x_{\infty}$ be any weak$^*$ limit point of the net $\tau_t(x)$ as $t 
\raro \infty$ which is an element in $\clg$ [Mo1]. In case (b) is true, we find that $x_{\infty} = \!E_{\clg}(x_{\infty}) = 
\phi_0(\!E_{\clg}(x))=\phi_0(x)1$. Thus $\phi_0(x)1$ is the unique limit point, hence weak$^*$ limit of $\tau_t(x)$ as 
$t \raro \infty$ is $\phi_0(x)1$. The equivalence statement for ergodicity also follows along the same line since the conditional 
expectation $\!E_{\cli}$ on the the von-Neumann algebra $\cli = \{ x: \tau_t(x)=x,\; t \ge 0 \}$ commutes with $(\tau_t)$ and thus 
satisfies $\!E_{\cli} \!E_{\clg} = \!E_{\clg} \!E_{\cli} = \!E_{\cli}$. This completes the proof that (a) and (b) are equivalent. 
That (b) and (c) are equivalent follows essentially along the same line since once more there exists a conditional expectation 
from $\clg$ to $\clg_0$ commuting with $(\tau_t)$ and any weak$^*$ limit point of the net $\tau_t(x)$ as $t$ diverges to infinity 
belongs to $\tau_s(\clg)$ for each $s \ge 0$, thus in $\clg_0$. We omit the details. \qed

\vsp
Now we investigate asymptotic behavior for quantum dynamical system dropping the assumption that $\phi_0$ is faithful. Let
$p$ be a sub-harmonic projection in $\cla_0$ for $(\tau_t)$ i.e.  $\tau_t(p) \ge p$ for all $t \ge 0$. Then $(\cla_0^p,
\tau^p_t,\phi^p_0)$ is a quantum dynamical semigroup where $\cla_0^p=p\cla_0p$ and $\tau^p_t(x)=p\tau_t(pxp)p$ for $x \in
\cla^p_0$ and $\phi^p_0(x)=\phi_0(pxp)$. In [Mo1] we have explored how ergodicity ( strong mixing ) of the original dynamics
can be determined by that of the reduced dynamics. Here we add one more result in that line of investigation.

\vsp
\NI {\bf THEOREM 2.4:} Let $(\cla_0,\tau_t,\phi_0)$ be a quantum dynamical systems with a normal invariant state $\phi_0$
and $p$ be a sub-harmonic projection for $(\tau_t)$. If $\mbox{s-limit}_{t \raro \infty}\tau_t(p)=1$
then the following statements are equivalent:

\NI (a) $||\phi\tau_t-\phi_0|| \raro 0$ as $t \raro \infty$ for any normal state on $\phi$ on $\cla_0$.

\NI (b) $||\phi^p \tau^p_t- \phi^p_0|| \raro 0$ as $t \raro \infty$ for any normal state $\phi^p$ on $\cla^p_0$.

\vsp
\NI {\bf PROOF: } That (a) implies (b) is trivial. For the converse we write
$||\phi \tau_t-\phi_0|| = \mbox{sup}_{x:||x|| \le 1} |\phi \tau_t(x)-\phi_0(x)|
\le \mbox{sup}_{\{x:||x|| \le 1\}} |\phi \tau_t(pxp)-\phi_0(pxp)| +
\mbox{sup}_{ \{x:||x|| \le 1 \} } |\phi \tau_t(pxp^{\perp})| +
\mbox{sup}_{ \{x:||x|| \le 1 \} } |\phi \tau_t(p^{\perp}xp)| +
\mbox{sup}_{ \{x:||x|| \le 1 \} } |\phi \tau_t(p^{\perp}xp^{\perp})|$.
Since $\tau_t((1-p)x) \raro 0$ in the weak$^*$ topology and
$|\phi \tau_t(xp^{\perp})|^2 \le |\phi \tau_t(xx^*)|\phi(\tau_t(p^{\perp}))| \le ||x||^2 \phi(\tau_t(p^{\perp})$
it is good enough if we verify that (a) is
equivalent to  $\mbox{sup}_{\{ x:||x|| \le 1 \} } |\phi\tau_t(pxp)-\phi_0(pxp)| \raro 0$
as $t \raro \infty$. To that end we first note that limsup$_{t \raro
\infty} \mbox{sup}_{x:||x|| \le 1} |\psi(\tau_{s+t}(pxp))-\phi_0(pxp)|$ is independent
of $s \ge 0$ we choose. On the
other hand we write $\tau_{s+t}(pxp) = \tau_s(p\tau_t(pxp)p) +
\tau_s(p\tau_t(pxp)p^{\perp})+ \tau_s(p^{\perp}\tau_t(pxp)p) +
\tau_s(p^{\perp}\tau_t(pxp)p^{\perp})$ and use the fact for any
normal state $\phi$ we have
$\mbox{limsup}_{t \raro \infty}\mbox{sup}_{x:||x|| \le 1}
|\psi(\tau_s(z\tau_t(pxp)p^{\perp})|\le ||z||\;|\psi(\tau_s(p^{\perp}))| $ for all $z \in \cla_0$.
Thus by our hypothesis on the support projection we conclude that (a) hold whenever (b) is true. \qed

\vsp
In case $\cla_0$ is a type-I von-Neumann algebra with center completely atomic, then $(\cla_0,\tau_t,\phi_0)$ 
is strong mixing if and only if $\clg = \!C$ [Mo1]. This criteria has been further explored in [Mo3] for an 
explicit necessary and sufficient condition on the coefficient associated with Stinespring representation [Da]
of a Markov map. This criteria in particular enable us to construct ergodic Markov map on a finite dimensional 
algebra with a faithful normal state and but not strong mixing. However in case the time variable is continuous and 
the von-Neumann algebra is the set of bounded linear operators on a finite dimensional Hilbert space $\clh_0$, 
by exploring Lindblad's representation [Li], Arveson [Ar] shows that a quantum dynamical semigroup with a faithful 
normal invariant state is ergodic if and only if the dynamics is strong mixing. In the following we prove a 
more general result exploring the criteria that we have obtained in Theorem 2.3. 

\vsp
\NI {\bf THEOREM 2.5:} Let $\cla_0$ be type-I with center completely atomic and $(\tau_t:t \in \!R)$ 
admits a normal state $\phi_0$. Then $(\cla_0,\tau_t,\phi_0)$ is strong mixing if and only if 
$(\cla_0,\tau_t,\phi_0)$ is ergodic. 

\vsp
\NI {\bf PROOF: } We first assume that $\phi_0$ is also faithful. We will verify now the criteria that 
$\clg_0$ is trivial when $(\tau_t)$ is ergodic. Since $\clg_0$ is remained invariant by the modular 
auto-morphism group associated with the faithful normal state $\phi_0$, by a theorem of Takesaki [Ta] there 
exists a faithful normal norm one projection from $\cla_0$ onto $\clg_0$. Now since $\cla_0$ is a von-Neumann 
algebra of type-I with center completely atomic, Stormer [So] says that $\clg_0$ is also type-I with center 
completely atomic. 

\vsp
Let $Q$ be a central projection in $\clg_0$. Since $\tau_t(Q)$ is also a projection and $\tau_t(Q) \raro Q$ as 
$t \raro 0$ we conclude that $\tau_t(Q)=Q$ for all $t \ge 0$ (center of $\clg$ being completely atomic and time 
variable $t$ is continuous ). Hence by ergodicity we conclude that $Q = 0$ or $1$. Hence $\clg_0$ can be identified 
with $\clb(\clk)$ for a separable Hilbert space $\clk$. Since $(\tau_t)$ on $\clb(\clk)$ is an automorphism we find a 
self-adjoint operator $H$ in $\clk$ so that $\tau_t(x)=e^{itH}xe^{-itH}$ for any $x \in \clb(\clk)$. Since it admits 
an ergodic faithful normal state, by [Fr, Mo1] we conclude that $\{ x \in \clb(\clk): x e^{itH} = e^{itH} x,\; t \in \!R \} 
= \IC$, which holds if and only if $\clk$ is one dimensional. Hence $\clg_0 = \IC$. 

\vsp
Now we deal with the general situation. Let $p$ be the support projection of $\phi_0$ in $\cla_0$.  So $p$ is 
a sub-harmonic projection in $\cla_0$ for $(\tau_t)$ i.e.  $\tau_t(p) \ge p$ for all $t \ge 0$. If $(\cla_0,\tau_t,\phi_0)$ 
is ergodic then by [Mo1,   ] $(\cla_p,\tau^p_t,\phi^p_0)$ is ergodic and $s-\mbox{limit}_{t \raro \infty}\tau_t(p)=1$ where 
$\cla_0^p=p\cla_0p$ and $\tau^p_t(x)=p\tau_t(pxp)p$ and $\phi^p_0(x)=\phi_0(pxp)$ for all $x \in \cla_0^p$. 
Since $\cla^p_0$ is also a type-I factor and $\phi_0^p$ is faithful, the reduced dynamics is strong mixing by 
the above argument. Now once more we appeal to [Mo1, Theorem  ] to complete the proof. \qed   

\vsp
We end this section with another application of Theorem 2.3, proving a result originated in [FNW1,FNW2,BJKW].

\vsp
\NI {\bf THEOREM 2.6: } Let $\cla_0$ be a type-I von-Neumann algebra with center completely atomic and $\tau$ be 
a completely positive map with a faithful normal invariant state $\phi_0$. Then the following are equivalent:

\NI (a) $(\cla_0,\tau_n,\phi_0)$ is strong mixing.

\NI (b) $(\cla_0,\tau_n,\phi_0)$ is ergodic and $\{ w \in S^1,\;\tau(x)=wx, \mbox{ for some non zero} x \in \cla_0 \} = 
\{ 1 \}$, where $S^1=\{w \in \IC: |w|=1 \}.$    

\vsp
\NI {\bf PROOF: } That `(a) implies (b)' is rather simple. To that end let $\tau(x)=wx$ for some $x \ne 0$ and $|w|=1$. 
Then $\tau^n(x)=w^nx$ and since the sequence $w^n$ has a limit point say $z, |z|=1$ we conclude by strong mixing that 
$zx=\phi_0(x)I$. Hence $x$ is a scaler and thus $x=\tau(x),\;x \ne 0$. So $w=1$. By taking $w=1$, we also get ergodic 
property, since by strong mixing $x=\phi_0(x)I$ for $\tau(x)=x$. 

\vsp
Now for the converse we will use our hypothesis that $\phi_0$ is faithful. To that end 
we plan to verify that $\clg_0$ is only scalers and appeal to Theorem 2.3 for strong mixing. Since there exists a conditional 
expectation from $\cla_0$ onto $\clg_0$, $\clg_0$ is once more a type-I von-Neumann algebra with center completely atomic. Let
$E$ be a non-zero atomic projection in the center of $\clg_0$. Since $\tau$ is an 
automorphism on $\clg_0$ each elements in the sequence $\{ \tau_k(E): k \ge 0 \}$ is an atomic projection in the center of 
$\clg_0$. If $\tau_n(E) \bigcap \tau_m(E) \ne 0$ and $n \ge m$ we find that $ \tau_m(\tau_{n-m}(E) \bigcap E) \neq 0$ 
and thus faithful and invariance property of $\phi_0$, we get $\phi(\tau_{n-m}(E) \bigcap E) > 0$. Once more by 
faithfulness we find $\tau_{n-m}(E) \bigcap E \ne 0$. So by atomic property of $E$ and $\tau_{n-m}(E)$ 
we conclude that $\tau_{n-m}(E)=E$. Thus either the elements in the infinite sequence $E, \tau(E),...., \tau^n(E)....$ 
are all mutually orthogonal or there exists an integer $n \ge 1$ so that the projections $E,\tau(E),.., \tau_{n-1}(E)$ 
are mutually orthogonal and $\tau^n(E)=E$. However for such an infinite sequence with mutually orthogonal projection 
we have $1 = \phi_0(I) \le \phi_0(\bigcup_{0 \le n \le m-1} \tau_n(E) ) = m \psi(E)$ for all $m \ge 1$. Hence $\psi(E)=0$ 
which is a contradiction, since $E$ is non-zero and $\phi_0$ is faithful.  

\vsp
Thus for any $w \in S^1$ with $w^n=1$, we have $\tau(x)=wx,$ where $x=\sum_{0 \le k \le n-1}w^k\tau_k(E) \ne 0 $. 
Hence by (b) we have $w=1$.  So $n=1$. In other words we have $\tau(E)=E$ for any atomic projection in the center of
$\clg_0$. Now by ergodicity we have $E=I$. Thus $\clg_0$ is a type-I factor isomorphic to $\clb(\clk)$ for some
Hilbert space $\clk$ and $\tau(x)=uxu^*$ for some unitary element in $\cla_0$. Since $(\clg_0,\tau_n,\phi_0)$ is 
ergodic by Theorem 2.3 we have $\{u,u* \}''=\clb(\clk)$, which holds if and only if $\clk$ is one dimensional ( check 
for an alternative proof that $\tau(u)=u$, thus $u=I$ by ergodicity and thus $\tau(x)=x$ for all $x \in \clg_0$ ).  
Hence $\clg_0=\IC$. This complete the proof that (b) implies (a). \qed

\newsection{ Minimal endomorphisms and Markov semigroups : }

\vsp
An E$_0$-semigroup $(\alpha_t)$ is a weak$^*$-continuous one-parameter semigroup
of unital $^*$-endomorphisms on a von-Neumann algebra $\cla$ acting on a Hilbert
space $\clh$. Following [Po1,Po2,Ar] we say $(\alpha_t)$ is {\it pure} if $\bigcap_{t \ge 0}
\alpha_t(\cla)=\IC$. For each $t \ge 0$, $\alpha_t$ being an endomorphism, $\alpha_t(\cla)$ is 
itself a von-Neumann algebra and thus $\bigcap_{t \ge 0}\alpha_t(\cla)$ is a limit of 
a sequence of decreasing von-Neumann algebras. Exploring this property Arveson proved that 
$(\alpha_t)$ is pure if and only if $||\psi_1\alpha_t-\psi_2\alpha_t|| \raro 0$ as $t \raro 
\infty$ for any two normal states $\psi_1,\psi_2$ on $\cla$. These criteria gets further simplified 
in case $(\alpha_t)$ admits a normal invariant state $\psi_0$ which says that $(\alpha_t)$ 
is pure if and only if $||\psi\alpha_t-\psi_0|| \raro 0$ as $ t \raro \infty$ for any normal 
state $\psi$. In such a case $\psi_0$ is the unique normal invariant state.  
However a pure $(\alpha_t)$ in general may not admit a normal invariant state [Po2,BJP] 
and this issue is itself an interesting problem.  

\vsp 
One natural question we wish to address here whether similar result is also true for a Markov 
semigroup $(\tau_t)$ defined on an arbitrary von-Neumann algebra $\cla_0$. This issue is already 
investigated in [Ar] where $\cla_0=\clb(\clh)$ and $(\tau_t)$ is assumed to be continuous in 
strong operator topology. He explored associated minimal dilation to an E$_0$-semigroups and thus 
make possible to prove that associated $E_0$-semigroup is pure if and only if 
$||\phi_1\tau_t-\phi_2\tau_t|| \raro 0$ as $t  \raro \infty$ for any two normal states 
$\phi_1,\phi_2$ on $\cla_0$. In case $(\tau_t)$ admits a normal invariant state the 
criteria gets simplified once more. In this section we will investigate this issue for an arbitrary 
von-Neumann algebra assuming that $(\tau_t)$ admits a normal invariant state $\phi_0$.

\vsp
To that end, we consider [Mo1] the minimal stationary weak Markov forward process 
$(\clh,F_{t]},j_t,\Omega,\; t \in \!R)$ and Markov shift $(S_t)$ associated with 
$(\cla_0,\tau_t,\phi_0)$ and set $\cla_{[t}$ to be the von-Neumann algebra generated 
by the family of operators $\{j_s(x): t \le s < \infty,\; x \in \cla_0 \}$. We recall 
that $j_{s+t}(x)=S^*_tj_s(x)S_t,\; t,s \in \!R$ and thus $\alpha_t(\cla_{[0}) \subseteq \cla_{[0}$ 
whenever $t \ge 0$. Hence $(\alpha_t,\; t \ge 0)$ is a E$_0$-semigroup on $\cla_{[0}$ with a 
invariant normal state $\Omega$ and 
$$ j_s(\tau_{t-s}(x))=F_{s]}\alpha_t(j_{t-s}(x))F_{s]} \eqno (3.1) $$ 
for all $x \in \cla_0$. We consider the GNS Hilbert space $(\clh_{ \pi_{\phi_0} }, \pi_{\phi_0}(\cla_0),\omega_0)$ 
associated with $(\cla_0,\phi_0)$ and define a Markov semigroup $(\tau_t^{\pi})$ on 
$\pi(\cla_0)$ by $\tau^{\pi}_t(\pi(x))= \pi(\tau_t(x)$. Furthermore we now identify $\clh_{\phi_0}$ as the 
subspace of $\clh$ by the prescription $\pi_{\phi_0}(x)\omega_0 \raro j_0(x)\Omega$. In such a case $\pi(x)$ 
is identified as $j_0(x)$ and aim to verify for any $ t \ge 0$ that
$$ \tau^{\pi}_t(PXP)=P\alpha_t(X)P \eqno (3.2) $$ for all $X \in \cla_{[0}$
where $P$ is the projection from $\clh$ on the GNS space. We use induction on $n \ge 1$. If $X=j_s(x)$ for some 
$s \ge 0$, (3.2) follows from (3.1). Now we assume that (3.2) is true for any element of the form $j_{s_1}(x_1)...j_{s_n}
(x_n)$ for any $s_1,s_2,...,s_n \ge 0$ and $x_i \in \cla_0$ for $1 \le i \le n$. Fix any $s_1,s_2,,s_n,s_{n+1} 
\ge 0$ and consider $X=j_{s_1}(x_1)...j_{s_{n+1}}(x_{n+1})$. Thus 
$P\alpha_t(X)P=j_0(1)j_{s_1+t}(x_1)...j_{s_{n+t}}(x_{n+1})j_0(1)$. If $s_{n+1} \ge s_n$,
we use (3.1) to conclude (3.2) by our induction hypothesis. Now suppose $ s_{n+1} \le s_n$. In that case if $s_{n-1} 
\le s_{n}$ we appeal to (3.1) and induction hypothesis to verify (3.2) for $X$. Thus we are left to consider 
the case where $s_{n+1} \le s_n \le s_{n-1}$ and by repeating this argument we are left to check only the case
where $s_{n+1} \le s_n \le s_{n-1} \le .. \le s_1 $. But $s_1 \ge 0=s_0$ thus we can appeal to (3.1) at the end of the
string and conclude that our claim is true for all elements in the $*-$ algebra generated by these elements 
of all order. Thus the result follows by von-Neumann density theorem.  We also 
note that $P=\tau^{\pi}_t(1)$ is a sub-harmonic projection [Mo1] for $(\alpha_t:t \ge 0)$ i.e. 
$\alpha_t(P) \ge P$ for all $t \ge 0$.  

\vsp
\NI {\bf THEOREM 3.1: } Let $(\cla_0,\tau_t,\phi_0)$ be a quantum dynamical semigroup with a normal invariant 
state for $(\tau_t)$. Then the GNS space $\clh_{\pi_{\phi_0}}$ associated with the normal state $\phi_0$ on $\cla_0 
$ can be realized as a closed subspace of a unique Hilbert space $\clh_{[0}$ up to isomorphism so that the following hold:

\NI (a) There exists a von-Neumann algebra $\cla_{[0}$ acting on $\clh_{[0}$ and a unital $*$-endomorphism 
$(\alpha_t,\; t \ge 0 )$ on $\cla_{[0}$ with a pure vector state $\phi(X)=<\Omega,X \Omega>$, $\Omega \in \clh_{[0}$
invariant for $(\alpha_t:t \ge 0)$. 

\NI (b) $P \cla P$ is isomorphic with $\pi(\cla_0)$ where $P$ is the projection onto $\clh_{\pi_{\phi_0}}$; 

\NI (c) $P\alpha_t(X)P=\tau^{\pi}_t(PXP)$ for all $t \ge 0$ and $X \in \cla_{[0}$; 

\NI (d) The closed span generated by the vectors $\{ \alpha_{t_n}(PX_nP)....\alpha_{t_1}(PX_1P)\Omega: 0 \le t_1 \le t_2 \le ..
\le t_k \le ....t_n, X_1,..,X_n \in \cla_{[0}, n \ge 1 \}$ is $\clh_{[0}$. 

\vsp
\NI {\bf PROOF:} The uniqueness up to isomorphism follows from the minimality property (d). \qed

\vsp
Following the literature [Vi,Sa,BhP,Bh] on dilation we say $(\cla_{[0},\alpha_t,\phi)$ is the minimal E$_0$semigroup associated with 
$(\cla_0,\tau_t,\phi_0)$. By a theorem [Ar, Proposition 1.1 ] we conclude that $\bigcap_{t \ge 0}\alpha_t(\cla_{[0}) = \IC$ if 
and only if for any normal state $\psi$ on $\cla_{[0}$, $||\psi\alpha_t-\psi_0|| \raro 0$ as $t \raro  \infty$, where 
$\psi_0(X)=<\Omega,X\Omega>$ for $X \in \cla_{0]}$. In the following proposition we explore that fact that $P$ is a sub-harmonic 
projection for $(\alpha_t)$ and by our construction $\alpha_t(P)=F_{t]} \uparrow 1$ as $t \raro \infty$. 

\vsp
\NI {\bf PROPOSITION 3.2:} $||\phi\tau^{\pi}_t-\phi_0|| \raro 0$ as $t \raro \infty$ for all normal state 
$\phi$ on $\pi(\cla_0)''$ if and only if $||\psi\alpha_t-\psi_0|| \raro 0$ as $t \raro \infty$ for all normal state 
$\psi$ on $\cla_{[0}$.

\vsp
\NI {\bf PROOF:} Since $F_{s]} \uparrow 1$ in strong operator topology by our construction and $\pi(\cla_0)$ is isomorphic to 
$F_{0]}\cla_{[0}F_{0]}$, we get the result by a simple application of Theorem 2.4. \qed   

\vsp
\NI {\bf THEOREM 3.3:} Let $\tau=(\tau_t,\;t \ge 0)$ be a weak$^*$ continuous Markov semigroup on $\cla_0$ 
with an invariant normal state $\phi_0$. Then there exists a weak$^*$ continuous E$_0$-semigroup 
$\alpha=(\alpha_t,\; t \ge 0)$ on a von-Neumann algebra $\cla_{[0}$ acting on a Hilbert space $\clh$ so that 
$$P\alpha_t(X)P=\tau^{\pi}_t(PXP),\; t \ge 0 $$ 
for all $X \in \cla_{[0}$, where $P$ is a sub-harmonic projection for $(\alpha_t)$ 
such that $\alpha_t(P) \uparrow I$.

\NI Moreover the following statements are equivalent:

\NI (a) $\bigcap_{t \ge 0}\alpha_t(\cla_{[0})=\!C$ 

\NI (b) $||\phi\tau_t^{\pi}-\phi_0|| \raro 0$ as $t \raro \infty$ for any normal state $\phi$ 
on $\pi(\cla_0)''$.

\NI (c) $\bigcap_{t \ge 0}\tau^{\pi}_t(\pi(\cla_0))=\!C$

\NI {\bf PROOF:} For convenience of notation we denote $\pi(\cla_0)''$ as $\cla_0$ in the following proof.
That (a) and (b) are equivalent follows by a Theorem of Arveson [Ar ] and Proposition 3.2.
Since $P\bigcap_{t \ge 0}\alpha_t(\cla_{[0})P=\bigcap_{t \ge 0}\tau^{\pi}_t(P\cla_{[0}P)=
\bigcap_{t \ge 0}\tau^{\pi}_t(\cla_0)$, if (a) is true then we have $\bigcap_{t \ge 0}\tau_t^{\pi}(\cla_0)
=\{zP: z \in \IC \}.$ Hence (c) is true. Conversely if (c) is true then $P\bigcap_{t \ge 0} \alpha_t(\cla_{[0})P
=\{zP: z \in \!C \}$ . Since $\bigcap_{t \ge 0}\alpha_t(\cla_{[0})$ is $(\alpha_t)$ invariant von-Neumann 
algebra, by homomorphism property we also get 
$\alpha_s(P) \bigcap_{t \ge 0}\alpha_t(\cla_{[0}) \alpha_s(P)= \{z\alpha_s(P):z \in \!C$. Since $\alpha_s(P) 
\uparrow 1$ as $s \raro \infty$ we conclude that (a) is also true. \qed

\vsp
Following [AM,Mo1] we say $(\clh,S_t,F_{t]},\Omega)$ is a {\it Kolmogorov's shift } if strong $\mbox{lim}_{t \raro - \infty}F_{t]}=|\Omega><\Omega|$. We
also recall here that Kolmogorov's shift property holds if and only if $\phi_0(\tau_t(x)\tau_t(y)) \raro  \phi_0(x)\phi_0(y)$ as $t \raro \infty$
for all $x,y \in \cla_0$. In such a case $\cla=\clb(\clh)$ [see the paragraph before Theorem 3.9 in [Mo1] ). If $\phi_0$ is faithful then $\cla_0$ 
and $\pi(\cla_0)$ are isomorphic, thus $\bigcap_{t \ge 0}\tau_t(\cla_0)=\!C$ if and only if $||\phi\tau_t-\phi_0|| \raro 0$ as $t \raro \infty$ for any normal state $\phi$ on $\cla_0$. 
Such a property is often called {\it strong ergodic property}. The following result says that there is a duality between strong ergodicity and
Kolmogorov's shift property.

\vsp
\NI {\bf THEOREM 3.4: } Let $(\cla_0,\tau_t,\phi_0)$ be a Markov semigroup with a faithful normal invariant state $\phi_0$. Then the
following are equivalent:

\NI (a) $\phi_0(\tilde{\tau}_t(x)\tilde{\tau}_t(y)) \raro \phi_0(x)\phi_0(y)$ as $t \raro \infty$ for any $x,y \in \cla_0$.

\NI (b) $||\phi \tau_t - \phi_0|| \raro 0$ as $t \raro \infty$ for any normal state $\phi$ on $\cla_0$.

\NI {\bf PROOF:} For each $t \in \!R$ let $\cla^b_{t]}$ be the von-Neumann algebra generated by the backward processes
$\{ j_s^b(x): - \infty <  s \le t \}$ [Mo1]. If (a) is true by Theorem 3.9 and Theorem 4.1 in [Mo1] we verify that weak$^*$ closure
of $\bigcup_{t \in \!R } \cla^b_{t]}$ is $\clb(\clh)$. Since for each $t \in \!R$ the commutant of $\cla^b_{t]}$ contains
$\cla_{[t}$ we conclude that $\bigcap_{t \in \!R} \cla_{t]}$ is trivial. Hence (b)
follows once we appeal to Theorem 3.3. For the converse, it is enough if we verify that $\phi_0(\tilde{\tau}_t(x)J\tilde{\tau}_t(y)J) \raro
\phi_0(x)\phi_0(y)$ as $t \raro \infty$ for any $x,y \in \cla_0$ with $y \ge 0$ and $\phi_0(y)=1$.  To that end we check the following
easy steps $\phi_0(\tilde{\tau}_t(x)J\tilde{\tau}_t(y)J) = \phi_0(\tau_t(\tilde{\tau}_t(x))JyJ)$ and for any normal state $\phi$,
$|\phi \circ \tau_t(\tilde{\tau}_t(x)) - \phi_0(x)| \le ||\phi \circ \tau_t - \phi_0|| ||\tilde{\tau}_t(x)|| \le ||\phi \circ \tau_t-\phi_0||
||x||$. Thus the result follows once we note that $\phi$ defined by $\phi(x)=\phi_0(xJyJ)$ is a normal state. \qed

\vsp
\NI {\bf THEOREM 3.5:}  Let $(\cla_0,\tau_t,\phi_0)$ be a Markov semigroup with a normal invariant state $\phi_0$. Consider 
the following statements:

\NI (a) $\phi_0(\tau_t(x)\tau_t(y)) \raro \phi_0(x)\phi_0(y)$ as $t \raro \infty$ for all $x,y \in \cla_0$.

\NI (b) the strong $\mbox{lim}_{t \raro - \infty}F_{t]} = |\Omega><\Omega|$.

\NI (c) $\cla = \clb(\clh)$ 

Then (a) and (b) are equivalent statements and in such a case (c) is also true. If $\phi_0$ is also faithful (c) is also 
equivalent to (a) ( and hence ( b)).  

\vsp
\NI {\bf PROOF: } That (a) and (b) are equivalent is nothing but a restatement of Theorem 3.9 in [Mo1]. That (b) implies (c) is
obvious since the projection $[\cla'\Omega]$, where $\cla'$ is the commutant of $\cla$, is the support of the vector state
in $\cla$. We will prove now (c) implies (a). In case $\cla=\clb(\clh)$, we have $\bigcap_{t \in \!R} \cla^b_{t]}= \!C$,
thus in particular $\bigcap_{t \le 0} \alpha_t(\cla^b_{0]})=\!C$. Hence by Theorem 3.3 applied for the time-reverse endomorphism 
we verify that $||\phi \tilde{\tau}_t-\phi_0|| \raro 0$ as $t \raro \infty$. Now (a) follows once we appeal to Theorem 3.4 
for the adjoint semigroups since $\tilde{\tilde{\tau}}_t=\tau_t$. \qed 

\vsp
Let $(\clb_0, \lambda_t,\;t \ge 0, \psi)$ be a unital $*-$endomorphism with an invariant normal state $\psi$ on a von-Neumann
algebra $\clb_0$ acting on a Hilbert space $\clk$. Let $P$ be the support projection for $\psi$. We set $\cla_0 = P \clb P$, a 
von-Neumann algebra acting on $\clh_0$, the closed subspace $P$, and $\tau_t(x)= P\lambda_t(PxP)P$, for any $x \in \cla_0$ and 
$t \ge 0$. Since $\lambda_t(P) \ge P$, it is simple to verify [Mo1] that $(\cla_0,\tau_t,\psi_0)$ is a
quantum dynamical semigroup with a faithful normal invariant state $\psi_0$, where $\psi_0(x)=\psi(PxP)$ 
for $x \in \cla_0$. Now we set $k_0(x)=PxP$ and $k_t(x)=\lambda_t(k_0(x))$ for $t \ge 0$. A routine verification 
says that $F_{s]}k_t(x)F_{s]}=k_s(\tau_{t-s}(x))$ for $0 \le s \le t$, where $F_{s]}=\lambda_s(P),\; s \ge 0$. Are these 
vectors $\{ \lambda_{t_n}(PX_nP)....\lambda_{t_1}(PX_1P)f: \; f \in \clh_0, 0 \le t_1 \le t_2 \le ..\le t_k \le ..t_n, X_1,..,
X_n \in \clb_0, n \ge 1 \}$ total in $\clh$? As an example we consider endomorphisms on $\clb(\clh)$ [BJP] 
with a pure mixing state, in such a case $\cla_0$ is only scalers thus the cyclic space associated with pure state is itself. 
Thus the problem is rather delicate even when the von-Neumann algebra is the algebra of all bounded operators on $\clk$. 
We will not address this problem here. Since $\lambda_t(P) \lambda_{t_n}(PX_nP)...\lambda_{t_1}(PXP) \clh_0 = 
\lambda_{t_n}(PX_nP)...\lambda_{t_1}(PXP) \Omega $ for $t \ge t_n$, $\mbox{lim}_{t \raro \infty}\lambda_t(P) =1$ is a 
necessary condition for cyclic property but not sufficient. However in the following we explore the fact the support 
projection $P$ is indeed an element in the von-Neumann algebra $\cln_0$ generated by the process 
$(k_t(x):\;t \ge 0,\;x \in \cla_0)$.

\vsp
To that end we consider little more general situation. Let $\clb_0$ be a $C^*$ algebra, $(\lambda_t:\;t \ge 0)$ be
a semigroup of endomorphisms and $\psi$ be an invariant state for $(\lambda_t:t \ge 0)$. We extend 
$(\lambda_t)$ to an automorphism on the $C^*$ algebra $\clb_{-\infty}$ of the inductive limit 
$$ \clb_0 \raro^{\lambda_t} \clb_0 \raro^{\lambda_t} \clb_0  $$
and extend also the state $\psi$ to $\clb_{-\infty}$ by requiring $(\lambda_t)$ invariance.
Thus there exists a directed set ( i.e. indexed by $\IT$ , by inclusion $\clb_{[-s} \subseteq \clb_{[-t}$ 
if and only if $t \ge s$ ) of C$^*$-subalgebras $\clb_{[t}$ of $\clb_{-\infty}$ so that the uniform closure of 
$\bigcup_{s \in \IT} \clb_{[s}$ is $\clb_{[-\infty}$. Moreover there exists an isomorphism
$i_0: \clb_0 \raro \clb_{[0}$ ( we refer [Sa] for general facts on inductive limit of C$^*$-algebras). 
It is simple to note that $i_t=\lambda_t \circ i_0$ is an isomorphism of $\clb_0$ onto $\clb_{[t}$ and 
$\psi_{-\infty} i_t =\psi$ on $\clb_0$. Let $(\clh_{\pi},\pi,\Omega)$ be the GNS space associated with 
$(\clb_{-\infty},\psi_{-\infty})$ and $(\lambda_t)$ be the unique normal extension to $\pi(\clb_{-\infty})''$. Thus 
the vector state $\psi_{\Omega}(X)=<\Omega, X \Omega>$ is $(\lambda_t)$ invariance and $(\pi(\clb_{[0})'',\lambda_t,
\;t \ge 0,\psi_{\Omega})$ is a quantum dynamics of endomorphism. Let $G_{0]}$ be the cyclic subspace of 
the vector $\Omega$ generated by $\pi(\clb_{[0})$. It is simple to check that $\pi(X)G_{0]}=G_{0]}\pi(X)G_{0]}$ 
for all $X \in \clb_{[0}$, hence each element in $\pi(\clb_{[0})''$ also commutes with $G_{0]}$. The the map 
$h: X \raro G_{0]}XG_{0]}$ is an homomorphism and the range is isomorphic to $\pi_0(\clb_0)''$, where $(\clh_0,\pi_0,\omega_0)$ 
be the GNS space associated with $(\clb_0,\psi)$. We identify the range of $h$ with $\pi_0(\clb_0)''$. It is simple
to verify that $h \circ \lambda_t(X)=\lambda_t(h(X))$ for all $X \in \pi(\clb_{[0})''$ and $t \ge $.

\vsp
Let $F_{t]}$ be the support projection of the normal vector state $\Omega$ in the von-Neumann sub-algebra 
$\pi(\clb_{[t})''$. $F_{t]}$ is a monotonically decreasing sequence of projections as $t \raro -\infty$. Let 
projection $Q$ be the limit. Thus $Q$ is the support projection for $\psi_{-\infty}$ in $\clb_{-\infty}$ and
$Q \ge |\Omega><\Omega|$. We aim to investigate when $Q$ is pure i.e. $Q=|\Omega><\Omega|$. 

\vsp
To that end we set von-Neumann algebra $\cln_0=F_{0]}\pi(\clb_{[0})''F_{0]}$ and 
define family $\{ k_t: \cln_0 \raro \pi(\clb_{-\infty})'',\; t \in \IT \}$ of $*-$homomorphisms by
$$k_t(x) = \lambda_t(F_{0]}xF_{0]}),\;\; x \in \cln_0$$  
It is a routine work to check that $(k_t:t \in \IT)$ is the unique up to isomorphism ( in the
cyclic space of the vector $\Omega$ generated by the von-Neumann algebra $\{k_t(x): t \in \IT, x \in \cln_0 
\}$ ) forward minimal weak Markov process associated with $(\cln_0,\eta_t,\psi_0)$ where 
$\eta_t(x)=F_{0]}\alpha_t(F_{0]}xF_{0]})F_{0]}$ for all $ t \ge 0$. Thus $Q=|\Omega><\Omega|$ when restricted 
to the cyclic space of the process if and only if $\psi_0(\eta_t(x)\eta_t(y)) \raro \psi_0(x)\psi_0(y) $ as 
$t \raro \infty$. In fact more is true. To that end let $P$ be the support projection of the vector state 
$\omega_0$ in von-Neumann algebra $\pi_0(\clb_0)''$ and $\cla_0=P\pi_0(\clb_0)''P$. We set $\tau_t(x)=
P\lambda_t(PxP)P$ for all $t \ge 0,\; x \in \cla_0$ and $\phi_0(x)=\psi(PxP)$.  

Thus $h(F_{0]})=P$ and by homomorphism property and commuting property with $(\lambda_t)$ we also check that 
$h(\cln_0)=\cla_0$ and $h(\eta_t(x))=h(F_{0]})\lambda_t(h(F_{0]})h(x)h(F_{0]}))= P \lambda_t(Ph(x)P)P=\tau_t(h(x))$
for all $t \ge 0$.   

\vsp
\NI {\bf THEOREM 3.6: } The following hold:

\NI (a) $k_t(I)=F_{t]}$ and $k_t(I)\pi(\clb_{[t})''k_t(I)= k_t(\cln_0)$ for all $t \in \IT$. 

\NI (b) $(\clh,k_t,F_{t]},\lambda_t,\Omega)$ is the minimal forward weak Markov process associated with 
$(\cln_0,\eta_t,\psi_0)$.  

\NI (c) $\psi_{-\infty}$ is a pure state if and only if $\phi_0(\tau_t(x)\tau_t(y)) \raro \phi_0(x)\psi_0(y)$ as $t \raro \infty$
for $x,y \in \cla_0$.

\vsp
\NI {\bf PROOF: } (a) is essentially by our construction and (b) is a routine work. We are left to prove only (c). For any 
fix $t \in \IT$ since $j_t(\cla_0)= F_{t]}\pi(\clb_{[t})''F_{t]}$,for any $X \in \clb_{[t}$ we have $QX\Omega=QF_{t]}XF_{t]}
\Omega=Qk_t(x)\Omega$ for some $x \in \cla_0$. Hence $Q=|\Omega><\Omega|$ if and only if $Q=|\Omega><\Omega|$ on the cyclic 
subspace generated by $ \{ k_t(x),\;t \in \IT, x \in \cla_0 \}$. Theorem 3.5 says now that $Q=|\Omega><\Omega|$ is and
only if $\psi_0(\eta_t(x)\eta_t(y)) \raro \psi_0(x)\psi_0(y)$ as $t \raro \infty$ for all $x \in \cln_0$, Since $h$ is an
homomorphism and $h \eta_t(x)= \tau_t(h(x))$, we also have $h(\eta_t(x))\eta_t(y))=\tau_t(h(x))\tau_t(h(x))$. Since 
$\phi_0 \circ h = \psi_0$ we conclude that result.  \qed

\newsection{ Sub-factors and Kolmogorov's shift: } 

\vsp
In this section we will investigate further the sequence of von-Neumann algebra $\cla_{[t}$ defined in the last section with an 
additional assumption that $\phi_0$ is also faithful.

\vsp
\NI {\bf THEOREM 4.1: } Let $(\cla_0,\tau_t,\phi_0)$ be a Markov semigroup with a faithful normal invariant state
$\phi_0$. If $\cla_0$ is a factor then $\cla_{[0}$ is a factor. Moreover 

\NI (a) $\cla_{[0}$ is a type-I (type-II, type-III) factor if and only if $\cla_0$ is a type-I 
(type -II , type-III) factor respectively.   

\NI (b) $\cla_0$ is a hyper-finite factor if and only if $\cla_{[0}$ is a hyper-finite factor.

\vsp
\NI {PROOF:} We first show factor property of $\cla_{[0}$. Note that the von-Neumann algebra $\cla^b_{0]}$ generated by the
backward process $\{ j^b_s(x):s \le 0, x \in \cla_0 \}$ is a sub-algebra of $\cla'_{[0}$, the commutant of $\cla_{[0}$. We fix 
any $X \in \cla_{[0} \bigcap \cla'_{[0}$ in the center. Then for any $y \in \cla_0$ we verify that 
$X j_0(y) \Omega= XF_{0]}j_0(y)\Omega=F_{0]}XF_{0]}j_0(y)\Omega = j_0(xy)\Omega$ for some $x \in \cla_0$. 
Since $Xj_0(y)=j_0(y)X$ we also have $j_0(xy)\Omega=j_0(yx)\Omega$. By faithfulness of the state $\phi_0$ we 
conclude $xy=yx$ thus $x$ must be a scaler. Thus we have $Xj_0(y)\Omega=c j_0(y)\Omega$ for some scaler $c \in \IC$. 
Now we use the property that $X$ commutes with forward process $j_t(x):\;x \in \cla_0, t \ge 0$ and as well as the 
backward processes $\{ j^b_t(x),\; t \le 0 \}$ to conclude that $X \lambda(t,x)= c \lambda(t,x)$. Hence $X=c$.

\vsp
Now if $\cla_0$ is a type-I factor, then there exists a non-zero minimal projection $p \in \cla_0$. In such a case we claim 
that $j_0(p)$ is also a minimal projection in $\cla_{[0}$. To that end let $X$ be any projection in $\cla_{[0}$ so
that $X \le j_0(p)$. Since $F_{0]}\cla_{[0}F_{0]}=j_0(\cla_0)$ we conclude that $F_{0]}XF_{0]}=j_0(x)$ for some $x \in \cla_0$. Hence
$X =j_0(p)Xj_0(p)=F_{0]}Xj_0(p)=j_0(xp)=j_0(px)$ 
Thus by faithfulness of the state $\phi_0$ we conclude that $px=xp$. Hence $X=j_0(q)$ where $q$ is a projection smaller then equal to $p$.
Since $p$ is a minimal projection in $\cla_0$, $q=p$ or $q=0$ i.e. $X=j_0(p)$ or $0$. So $j_0(p)$ is also a minimal projection.
Hence $\cla_{[0}$ is a type-I factor. For the converse statement we trace the argument in the reverse direction. Let $p$ be a 
non-zero projection in $\cla_0$ and claim that there exists a minimal projection $q \in \cla_0$ so that $0 < q \le p$. Now 
since $j_0(p)$ is a non-zero projection in a type-I factor $\cla_{[0}$ there exists a non-zero projection $X$ which is minimal in $\cla_{[0}$
so that $0 < X \le j_0(p)$. Now we repeat the argument to conclude that $X=j_0(q)$ for some projection $q$. Since $X \neq 0$ and minimal, 
$q \neq 0$ and minimal in $\cla_0$. This completes the proof for type-I case. We will prove now the case for Type-II. 

\vsp
Let $\cla_{[0}$ be type-II then there exists a finite projection $X \le F_{0]}.$ Once more $X=F_{0]}XF_{0]}=j_0(x)$ for some projection 
$x \in \cla_0$. We claim that $x$ is finite. To that end let $q$ be another projection so that $q \le x$ and $q=uu^*$ and $u^*u=x$. Then $j_0(q) 
\le j_0(x)=X$ and $j_0(q)=j_0(u)j_0(u)^*$ and $j_0(x)=j_0(u)^*j_0(u)$. Since $X$ is finite in $\cla_{[0}$ we conclude that $j_0(q)=j_0(x)$. By 
faithfulness of $\phi_0$ we conclude that $q=x$, hence $x$ is a finite projection. Since $\cla_0$ is not type-I, it is type-II.
For the converse let $\cla_0$ be type-II. So $\cla_{[0}$ is either type-II or type-III. We will rule out that the possibility 
for type-III. Suppose not, i.e. if $\cla_{[0}$ is type-III, for every projection $p \ne 0$, there exists $u \in \cla_{[0}$ so
that $j_0(p)=uu^*$ and $F_{0]}=u^*u$. In such a case $j_0(p)u=uF_{0]}$. Set $j_0(v)=F_{0]}uF_{0]}$ for some $v \in \cla_0$. 
Thus $j_0(pv)=j_0(v)$.  Once more by faithfulness of the normal state $\phi_0$, we conclude $pv=v$. So $j_0(v)=uF_{0]}$. 
Hence $j_0(v^*v)=F_{0]}$. Hence $v^*v=1$ by faithfulness of $\phi_0$. Since this is true for any non-zero projection $p$ in
$\cla_0$, $\cla_0$ is type-III, which is a contradiction. Now we are left to show the statement for type-III, which is true
since any factor needs to be either of these three types. This completes the proof for (a).

\vsp
For (b) we recall for a factor, hyperfinite is equivalent to being generated by an ascending sequence 
of finite dimensional von-Neumann algebras [BR,El]. Let $\cla_0$ be hyperfinite and $\{ \cln^n:n \ge 1 \}$ 
be such a sequence of finite dimensional von-Neumann algebras. For each $n \ge 1$ we set von-Neumann
sub-algebras $\cln^n_{[0} \subseteq \cla_{[0}$ generated by the elements $\{ j_t(\cln^n): t = {r \over {2^n}}, 0 \le r 
\le n2^n \}$. Any arbitrary product of elements in the set is reduced to a product of elements of at most 
$n2^n$ elements from the set $\{ j_t(\tau_s(\cln^n)): \mbox{where }\; s,\;t \in \{ {r \over {2^n}}, 0 \le r
\le n2^n \} \}$. Thus each $\cln^n_{[0}$ is finite dimensional and ascending with $n$. By the weak$^*$ continuity of 
Markov semigroup we check that the sequence generates $\cla_{[0}$. Hence by our earlier remark $\cla_{[0}$ is hyperfinite. 
For the converse we recall for a factor $\clm$ acting on a Hilbert space $\clh$, Tomiyama's property 
( i.e. there exists a norm one projection $E: \clb(\clh) \raro \clm$, see [BR1] page-151 for details ) 
is equivalent to hyperfinite property. For a hyperfinite factor $\cla_{[0}$, $j_0(\cla_0)$ is a factor in the GNS space
identified with the subspace $F_{0]}$. Let $E$ be the norm one projection from $\clb(\clh_{[0})$ on $\cla_{[0}$ and 
verify that the completely positive map $E_0: \clb(\clh_0) \raro \cla_0$ defined by $E_0(X)= F_{0]}E(F_{0]}XF_{0]})F_{0]}$ 
is a norm one projection from $\clb(F_{0]})$ to $\cla_0$. This completes the proof. \qed

\vsp
Let $H$ be a Hilbert space, $B(H)$ the algebra of bounded operators on $H$, and $\cle$ a complete Boolean algebra (complete orthocomplemented 
distributive lattice) with minimal element 0 and maximal element 1. $I \neq 0$ in $\cle$ is an atom if $J < I$ implies that $J=0$. $\cle$ is atomic if 
for every $J \in \cle$ there is an atom $I \leq J$; $\cle$ is continuous if it has no atom. A complete Boolean algebra of factors is a mapping 
$I \rightarrow R(I)$ from $\cle$ into the von Neumann algebras on $H$, such that $R(I')=R(I)'$, $R(\bigwedge I_\alpha)=\bigcap
R(I_\alpha)$, $R(\bigvee I_\alpha)=(\bigcup R(I_\alpha))''$, $R(1)=B(H)$, and, for every $I \in \cle$, $R(I)$ is a factor ($I'$ denotes the 
complement of $I$, and $R(I)'$ the commutant of $R(I)$).

\vsp
We set family $\cla_{[s,t)}=\cla_{[s} \bigcap \cla'_{[t},\; -\infty < s \le t < \infty $ of factors, $\cla_{[s,\infty)}=\cla_{[s}$ and 
$\cla_{(-\infty,t]}=\cla'_{[t}.$ 

\vsp
\NI {\bf THEOREM 4.2:} The map $[s,t) \raro \cla_{[s,t)}$ has a unique extension to a complete boolean algebra of factors if and only if 
$F_{-t]} \raro |\Omega><\Omega|$ and $F_{[t} \raro |\Omega><\Omega|$ as $t \raro \infty$.     

\vsp
\NI{\bf Proof:} By Theorem 3.5 we have $F_{-t]} \raro |\Omega><\Omega|$ as $t \raro \infty$ if and only if $(\bigcup_{t \in \!R}\cla_{[t})''
=\clb(\clh)$. By duality we also have $F_{[t} \raro |\Omega><\Omega|$ as $t \raro \infty$ if and only if $\bigcap_{t \in \!R} \cla_{[t}=
\IC$. Thus $\cla_{[s,t)} \uparrow \clb(\clh)$ as $[s,t) \uparrow (-\infty,\infty)$ if and only if $F_{-t]} \raro |\Omega><\Omega|$ and
$F_{[t} \raro |\Omega><\Omega|$ as $t \raro \infty$. \qed

\newsection{ Complete boolean algebra of type-I factors: }

\vsp
For a type-I factor $\cla_0$, $\cla_{[0}$ is also a type-I factor, thus there exist Hilbert spaces $\clf_{0)}$ and $\clf_{[0}$ so that 
$\clh$ is isomorphic to $\clf_{0)} \otimes \clf_{[0}$ and $\cla_{[0}$ is isomorphic to $I \otimes \clb(\clf_{[0})$. Since for any $t \ge 0$, 
$\cla_{[t}$ is a type-I sub-factor of $\cla_{[0}$, thus via isomorphism is also a type-I sub-factor of $\clb(\clf_{[0})$. Thus we also 
find a Hilbert space $\clf_{[0,t)}$ so that $\clf_{[0}$ is isomorphic with $\clf_{[0,t)} \otimes \clf_{[t}$ and $\cla_{[t}$ is 
isomorphic with $I \otimes \clb(\clf_{[t}).$ Moreover $\clf_{[0,s)} \otimes \clf_{[s,t)}$ is isomorphic with $\clf_{[0,t)}$ for any $s < t$. Since $\cla_{[0}$ 
is isomorphic to $\cla_{[s}$, we also verify that $\clf_{[s,t)}$ is isomorphic to $\clf_{[0,t-s)}$. Thus the family $\clp_t=\clf_{[0,t)},\; t 
\ge 0$ is a product system [Ar] of Hilbert spaces in $\clf_{[0}$ i.e. $P_t \otimes P_s$ is isomorphic to $P_{s+t}$ for any $s,t \ge 0$.

\vsp
Moreover for a type-I factor $\cla_0$, Kolmogorov's property of a Markov semigroup $(\cla_0,\tau_t,\phi_0)$ is equivalent to strong 
mixing. Since strong mixing property is time reversible, by duality Kolmogorov property of the adjoint Markov semigroup is also 
equivalent to strong mixing. Thus by Theorem 4.2 we conclude that the map $[s,t) \raro \cla_{[s,t)}$ has a unique extension to a 
complete boolean algebra of type-I factors if and only if $(\cla_0,\tau_t,\phi_0)$ is strongly mixing. In such a case the pure vector 
state $\phi$ on $\clb(\clh)$ is quasi-equivalent [BR] to product states $\phi_- \otimes \phi_+$ where $\phi_+$ and $\phi_{-}$ are the 
normal states $\phi$ restricted to $\cla_{[0}$ and $\cla'_{[0}$ ( the commutant ) respectively.

\vsp
\NI {\bf THEOREM 5.1: } Let $(\cla_0,\tau_t,\phi_0)$ be as in Proposition 4.1 and $\cla_0$ be a type-I factor. Then the following hold:

\NI (a) There exists complex separable Hilbert spaces $\clf_{0)}, \clf_{[0}$ and an unitary operator $U_0: \clf_{0)} \otimes \clf_{[0} \raro \clh$ 
so that
$U_0^*\cla_{[0}U_0= I_{\clf_{0)}} \otimes \clb(\clf_{[0})$,

\NI (b) $(\clb(\clf_{[0}),\beta_t,\;t \ge 0,\psi)$ is an unital $*$-endomorphisms, where $\beta_t(X) = U_0^*\alpha_t(U_0 (I_{\clf_{0)}} \otimes
X)U_0^*)U_0,\;\; t \ge 0$ and $X \in \clb(\clf_{[0})$ and the normal state $\psi(X) = \phi(U_0(I_{\clf_{0)}} \otimes X)U_0^*),\;\;
X \in \clb(\clf_{[0})$ is invariant for $(\beta_t);$ 

\NI (c) Let $P_0$ be the support projection in $\clf_{[0}$ of $\psi$ and $\clk_0$ be the Hilbert subspace $P_0$, then 

\vsp
\NI (i) $P_0\clb(\clf_{[0})P_0$ is isomorphic to $\pi_{\phi_0}(\cla_0)$, where $\pi_{\phi_0}$ is the GNS representation associated with $\phi_0;$

\NI (ii) $U^*_0j_0(x)U_0 \equiv I_{\clf_{0)}} \otimes \pi_{\phi_0}(x)$;

\NI (iii) $\pi_{\phi_0}(\tau_t(x)) \equiv P_0\beta_t(P_0 \pi_{\phi_0}(x)P_0)P_0$ for any $x \in \cla_0$ and $t \ge 0$;  

\NI (d) The von-Neumann algebra generated by $\{\beta_t(P_0xP_0):t \ge 0, x \in \clb(\clh_0) \}$ is $\clb(\clf_{[0})$;

\NI (e) The set $\{ \beta_{t_n}(P_0x_nP_0)...\beta_{t_1}(P_0x_1P_0)f: f \in \clh_0,x_1,,,x_n \in \clb(\clh_0),0 \le t_1 \le ...\le t_n,\; n \ge 1 \}$
is total in $\clf_{[0}$; 

\NI (f) $U^*_0F_{[0} \clb(\clh) F_{[0}U_0 \equiv \pi_{\phi_0}(\cla_0) \otimes \clb(\clf_{[0})$;

\vsp
\NI{\bf PROOF:} (a) follows since $\cla_{[0}$ is also a type-I factor by Proposition 4.1. (b) is simple to verify. For (c) we recall
$j_0(1)=F_{0]} \in \cla_{[0}$. Hence $U^*_0F_{0]}U_0=I_{\clf_{(0}} \otimes P$ where $P$ is a projection in $\clf_{0]}$. However
$\psi(P)=1$, hence $P \ge P_0$. We claim that $P=P_0$. To that end note that $U_0(I_{\clf_{(0}} \otimes P_0)U^*_0 \le F_{0]}$, hence 
$U_0(I_{\clh_{(0}} \otimes P_0)U^*_0=j_0(x)$ for some projection $x \in \cla_0$. But $\phi(x)=\phi_0(x^*x)=1$, hence $x=1$ by faithfulness 
of $\phi_0$. Thus $U^*_0F_{0]}U_0=I_{\clf_{(0}} \otimes P_0$ and $U^*_0j_0(x)U_0 = P_0xP_0$, where we have identified $\cla_0$ with $\clb(\clk_0)$.
Now it is routine to verify (c) using (b). (d) is rather obvious now by (a) and statement (i) of (c). (e) is trivial once we use (d). 
(f) is rather delicate. To that end first note that $F_{[0} \in \cla'_{[0}$, the commutant of $\cla_{[0}$, in fact little more is
true, $F_{[0} \in \cla^b_{0]}$, the von-Neumann algebra generated by the the backward processes $\{ j^b_t(x): t \le 0,x \in \cla_0 \}$
and $\cla^b_{0]} \subseteq \cla'_{[0}$. Thus $U^*_0F_{[0}U_0=Q_0 \otimes I_{\clf_{[0}}$ where $Q_0$ is a projection in $\clf_{(0}$. Now 
we follow the steps in the proof of (c) applied to the backward process to conclude that $Q_0$ is the support projection for $\psi$
once restricted to $\cla^b_{0]}$. It is simple to note that $U^*_0\cla^b_{0]}U_0 \subseteq \clb(\clf_{(0}) \otimes I_{\clf_{0]}}$. We 
claim that the equality hold. This follows once we recall from the proof of Proposition 4.1 that von-Neumann algebra $\cla_{[0}$
together with $\cla^b_{0]}$ generate $\clb(\clh)$. Since $Q_0\clb(\clf_{(0})Q_0$ is isomorphic to $\pi_{\phi_0}(\cla_0))$, (f) 
follows. \qed

\newsection{ Pure state on the two sided quantum spin chain: }

\vsp
In this section we essentially recall basic facts on Cuntz algebras presented as in [BJKW] and 
investigate when a translation invariant state on quantum spin chain is pure. This in particular answers 
an important question how Kolmogorov's property of the associated Popescu system is related with purity 
of the state. Perhaps it is the most interesting application of our result obtained in section 3. 
 
\vsp
First we recall that if $d \in \{2,3,.., \}$, the Cuntz algebra is the universal $C^*$-algebra generated by elements $\{s_1,s_2,...,s_d \}$ subject
to the relations:

$$s^*_is_j=\delta^i_j1$$
$$\sum_{1 \le i \le d } s_is^*_i=1.$$

\vsp
There is a canonical action of the group $U(d)$ of unitary $d \times d$ matrices on $\clo_d$ given by
$$\beta_g(s_i)=\sum_{1 \le j \le d}\overline{g^j_i}s_j$$
for $g=((g^i_j) \in U(d)$. In particular the gauge action is defined by
$$\beta_z(s_i)=zs_i,\;\;z \in \IT \subset \IC $$
If UHF$_d$ is the fixed point subalgebra under the gauge action, then UHF$_d$ is the closure of the
linear span of all wick ordered monomials of the form
$$s_{i_1}...s_{i_k}s^*_{j_k}...s^*_{j_1}$$
which is also isomorphic to the UHF$_d$ algebra
$$M_{d^\infty}=\otimes^{\infty}_1M_d$$
so that the isomorphism carries the wick ordered monomial above into the matrix element
$$e^{i_1}_{j_1}(1)\otimes e^{i_2}_{j_2}(2) \otimes....\otimes e^{i_k}_{j_k}(k) \otimes 1 \otimes 1 ....$$
and the restriction of $\beta_g$ to $UHF_d$ is then carried into action
$$Ad(g)\otimes Ad(g) \otimes Ad(g) \otimes ....$$

\vsp
We also define the canonical endomorphism $\lambda$ on $\clo_d$ by
$$\lambda(x)=\sum_{1 \le i \le d}s_ixs^*_i$$
and the isomorphism carries $\lambda$ restricted to UHF$_d$ in the one-sided shift
$$y_1 \otimes y_2 \otimes ... \raro 1 \otimes y_1 \otimes y_2 ....$$
on $\otimes^{\infty}_1 M_d$. Note that $\lambda \beta_g = \beta_g \lambda $ on UHF$_d$.

\vsp
Let $d \in \{2,3,..,,..\}$ and $\IZ_d$ be a set of $d$ elements.  $\cli$ be the set of finite sequences $I=(i_1,i_2,...,i_m)$ where
$i_k \in \IZ_d$ and $m \ge 1$. We also include empty set $\emptyset \in \cli$ and set $s_{\emptyset }=1=s^*_{\emptyset}$,
$s_{I}=s_{i_1}......s_{i_m} \in \clo_d $ and $s^*_{I}=s^*_{i_m}...s^*_{i_1} \in \clo_d$. In the following we recall from [BJKW]
a crucial result originated in [Po,BJP].

\vsp
\NI{\bf THEOREM 6.1: } There exists a canonical one-one correspondence between the following objects:

\NI (a) States $\hat{\omega}$ on $\clo_d$

\NI (b) Function $C: \cli \times \cli \raro \IC$ with the following properties:

\NI (i) $C(\emptyset, \emptyset)=1$;
\NI (ii) for any function $\lambda:\cli \raro \IC$ with finite support we have
         $$\sum_{I,J \in \cli} \overline{\lambda(I)}C(I,J)\lambda(J) \ge 0$$

\NI (iii) $\sum_{i \in \IZ_d} C(Ii,Ji) =C(I,J)$ for all $I,J \in \cli$.

\NI (c) Unitary equivalence class of objects $(\clk,\Omega,V_1,..,V_d)$ where

\NI (i) $\clk$ is a Hilbert space and $\Omega$ is an unit vector in $\clk$;
\NI (ii) $V_1,..,V_d \in \clb(\clk)$ so that $\sum_{i \in \IZ_d} V_iV^*_i=1$;
\NI (iii) the linear span of the vectors of the form $V^*_I\Omega$, where $I \in \cli$, is dense in $\clk$.

\vsp
Where the correspondence is given by a unique completely positive map $R: \clo_d \raro \clb(\clk)$ so that

\NI (i) $R(s_Is^*_J)=V_IV^*_J;$

\NI (ii) $\hat{\omega}(x)=<\Omega,R(x)\Omega>;$

\NI (iii) $\hat{\omega}(s_Is^*_J)=C(I,J)=<V^*_I \Omega,V^*_J\Omega>.$

\NI (i) For any fix $g \in U_d$ and the completely positive map $R_g:\clo_d \raro \clb(\clk)$ defined by $R_g= R
\circ \beta_g$ give rises to a Popescu system given by $(\clk,\Omega,\beta_g(V_i),..,\beta_g(V_d))$ where
$\beta_g(V_i)=\sum_{1 \le j \le d} \overline{g^i_j} V_j.$

\vsp
Let $\omega$ be a translation invariant ergodic states (extremal states) on UHF$_d$ algebra $\otimes_{\IN} M_d$. Following
[BJKW, section 7], we consider the set
$$K_{\omega}= \{ \psi: \psi \mbox{ is a state on } \clo_d \mbox{ such that } \psi \lambda = \psi \mbox{ and } \psi_{|\mbox{UHF}_d} 
= \omega \}$$ $K_{\omega}$ is a non empty convex and compact in weak topology. $K_{\omega}$ is a face in the $\lambda$ invariant 
states since $\omega$ is extremal. We recall Lemma 7.4 of [BJKW] in the following proposition.

\vsp
\NI {\bf PROPOSITION 6.2:} $\psi \in K$ is an extremal points in $K$ if and only if $\psi$ is a factor state and moreover
all other extremal points have the form $\psi \beta_z$ for some $z \in \IT$.

\vsp
\NI {\bf PROPOSITION 6.3:} Let $\omega$ be an extremal point in UHF$_d$ C$^*$ algebra and then there exists a
von-Neumann algebra $\clm$ acting on a Hilbert space $\clk$ so that the following hold:

\NI (a) $\clm$ is a factor.

\NI (b) There exists $V_1,V_2,...,V_d$ bounded operators on $\clm$ so that $\tau^v(X)=\sum_{1 \le i \le d} V_iXV^*_i$ on
$\clm$ is an ergodic map with a faithful normal invariant state $\phi_0$ on $\clm$.

\NI (c) For any $I=(i_1,i_2,...,i_k),J=(j_1,j_2,...,j_k)$ with $|I|=|J| < \infty$ we have
$\omega(e^{i_1}_{j_1} \otimes e^{i_2}_{j_2} \otimes e^{i_3}_{j_3} \otimes ...\otimes e^{i_k}_{j_k}) = \phi_0(V_IV^*_J)$.

Conversely any $\lambda$ invariant state $\omega$ on UHF$_d$ satisfying (a)-(c) is an extremal point. 

\vsp
\NI {\bf PROOF:} We fix an extremal point $\psi$ in $K_\omega$ and consider the GNS space  $(\clh_{\pi_\psi},
\pi_{\psi}(\clo_d),\Omega_{\psi})$ associated with $(\clo_d,\psi)$. Set $S_i=\pi_{\psi}(s_i)$ and consider
the normal state $\psi_{\Omega}$ on $\pi(\clo_d)''$ defined by $\psi_{\Omega}(X) = <\Omega, X\Omega>$. Thus
$(\pi(\clo_d)'',\Lambda,\psi)$ is an ergodic Markov map where
$$\Lambda(X)=\sum_{1 \le i \le d}S_iXS^*_i. $$
Let $P$ be the support projection in $\pi_{\psi}(\clo_d)''$ of the normal state $\psi_{\Omega}$. Thus $P$ be a sub-harmonic
projection for $\beta$, thus by [FR,Mo1] we have
\be
PS^*_iP=S^*_iP
\ee
and the reduce dynamic
\be
\tau^v(x)=P\alpha(PxP)P
\ee
on $\clm=P \pi_{\psi}({\clo_d})''P$ is also ergodic with a faithful normal state $\phi$ where $\phi_0(x)=\psi_{\Omega}(PxP)$ for all $x \in \clm$.
Hence $\clm=\{V_i,V^*_i \}''$ and $\clm$ is a factor, where we defined $V_i=PS_iP$ for all $1 \le i \le d$. That (c) is satisfied 
follows by the relation (6.1).
Conversely for any given family of Popescu system $(\clk,\clm,V_i,\phi_0)$ we consider the minimal dilation  $(\clh,S_i,P)$ as described
in [BJKW, Theorem 5.1] where

\NI (a) $\clh$ is a Hilbert space and $P$ is a projection so that $P\clh=\clk$;

\NI (b) The family of isometric operators $(S_i,\;1 \le i \le d)$ satisfies Cuntz's relation $\sum_{1 \le i \le d}S_iS_i^*=1$;

\NI (c) $(1-P)S_i^*P=0$ and $P$ is cyclic for the representation i.e. the set $\{ S_If: |I| < \infty, Pf=f,\;f \in 
\clh \}$ is total in $\clh$.

\vsp
We define a unique state $\hat{\omega}$ on $\clo_d$ by prescribing $\hat{\omega}(s_Is^*_J)=\phi(V_IV^*_J)$. Thus the GNS representation
associated with $(\clo_d,\hat{\omega})$ identifies with $\pi_{\hat{\omega}}(s_i) = S_i$. Now (c) guarantees
that $\Lambda_n(P) \uparrow I$ as $n \uparrow \infty$. Hence Theorem 3.6 in [Mo1] ensures that the endomorphism 
$(\pi_{\hat{\omega}}(\clo_d)'', \Lambda_n,\phi_{\omega})$ is ergodic since the reduced dynamics $(\clm,\tau^v_n,\phi)$ 
is so. Thus $\hat{\omega}$ is an extremal $\lambda$ invariant state on $\clo_d$.  \qed

\vsp
Let $\omega'$ be a translation invariant ergodic state on UHF$_d$ algebra $\otimes_{\IZ}M_d$ and $\omega$ be the 
restriction of $\omega'$ to UHF$_d$ algebra $\clb_0=\otimes_{\IN} M_d$. We fix any $\psi \in K_{\omega}$ an extremal point 
and consider the associated Popescu system as described in Proposition 6.3. Then a simple application of Theorem 3.6 
says that the inductive limit state $\psi_{-\infty}$ on the inductive limit $(\clo_d,\psi) \raro^{\lambda} (\clo_d,\psi) 
\raro^{\lambda} (\clo_d,\psi)$ is pure if and only if $\phi_0(\tau^v_n(x)\tau^v_n(y)) \raro \phi_0(x)\phi_0(y)$ for 
all $x,y \in \clm$ as $n \raro \infty$. 

\vsp
The von-Neumann algebra $\{ S_IS^*_J: |I|=|J| < \infty \}''$ acts on the cyclic subspace 
$\clh_{\pi_0}$ generated by the vector $\Omega$. This is isomorphic with the GNS representation of associated 
with $(\clb_0,\omega)$. The inductive limit $(\clb_{-\infty},\phi_{-\infty})$ described 
as in Proposition 3.6 associated with $(\clb_0,\lambda_n,\;n \ge 0,\omega)$ is UHF$_d$ algebra 
$\otimes_{\IZ}M_d$ and the inductive limit state is $\omega'$. Let $Q$ be the support projection of the state 
$\psi$ in $\pi_0(\clb_0)''$ and $\cla_0=Q\pi(\clb_0)''Q$. Since $\Lambda(Q) \ge Q$, from the identity $Q\Lambda(I-Q)Q=0$ we 
conclude that $(I-Q)S_i^*Q=0$ for all $1 \le i \le d$. However $Q \ge P$, where $P$ is the support projection of $\psi$ on 
$\pi(\clo_d)''$. We set $l_i=QS_iQ$ and note that $\psi(S_IS^*_J)=\psi_0(l_Il^*_J)$ for all $|I|=|J| < \infty$, where 
$\psi_0(x)=\psi(QxQ)$ for $x \in \cla_0$. Since $\bigcap_{n \ge 1}\Lambda_n(\pi_0(\clb_0)'')=\IC$, we have 
$||\Lambda^n(X)-\phi(X)I|| \raro 0$ as $n \raro \infty$, in particular, $\Lambda_n(Q) \uparrow I$ as $n \raro \infty$. 
Hence $ \{ S_If: |I| < \infty,\; Qf=f,f \in \clh_{\pi} \}$ is total in $\clh_{\pi_0}$.

\vsp
\NI {\bf THEOREM 6.4: } Let $\omega'$ be a translation invariant ergodic state on UHF$_d$ algebra $\otimes_{\IZ}M_d$. Then 
the following are equivalent:

\NI (a) $\omega'$ is a pure state.

\NI (b) $\psi_0(\eta_n(x)\eta_n(y)) \raro \psi_0(x)\psi_0(y)$ as $n \raro \infty$ for all $x,y \in \cla_0$, where
$\eta(x)=\sum_il_ixl_i^*$ for all $x \in \cla_0$.

\NI {\bf PROOF:} We consider the dynamics $(\clb_0,\lambda_n,\psi)$ and appeal to Theorem 3.6. \qed   

\vsp
That for a given $\omega'$, the Popescu system is uniquely determined modulo a unitary equivalence follows 
from Proposition 6.2 and Proposition 6.3. Thus the Markov semigroup $(\cla_0,\eta,\psi_0)$ is also uniquely
determined modulo a unitary conjugation. Thus the criterion appeared in Theorem 6.4 (b) is independent of the 
extremal point $\psi \in K_{\omega}$ that we have chosen.

\vsp
\NI {\bf COROLLARY 6.5: } Let $\omega'$ be as in Theorem 6.4 and its restriction $\omega$ to UHF$_d$ algebra
$\otimes_{\IN}M_d$ be a type-I factor state. Then $\omega'$ is a pure state if and only if $(\cla_0,\eta_n,\psi_0)$ 
is strongly mixing. In such a case $\cla_0=\clm$ and $l_k=V_k$ for all $1 \le k \le d$.  

\NI {\bf PROOF: } Since strong mixing is equivalent to Kolmogorov's property for a type-I von-Neumann algebra 
( see Theorem 4.7 in [Mo1]) the first part of the corollary follows from Theorem 6.4. 

\vsp
Since $\pi_{\omega}(\otimes_{\IN}M_d$ is a type-I factor, the unique canonical normal endomorphism 
$\Lambda:\pi_{\omega}(\otimes_{\IN}M_d))'' \raro \pi_{\omega}(\otimes_{\IN}M_d)''$ has Powers index $d$ 
and $\Lambda(X)= \sum_{1 \le k \le d}S_iXS^*_i$ for all $X \in \pi_{\omega}(\otimes_{\IN}M_d)''$
where $(S_i)$ are elements in $\pi_{\omega}(\otimes_{\IN}M_d)''$ satisfying Cunz's relations. 
For more details we refer to [BJP]. Since $\omega$ is an ergodic state, $\Lambda_n(Q) \uparrow I$ as $n \uparrow I$, 
where $Q$ is the support projection of $\omega$ in $\pi_{\omega}(\otimes_{\IN}M_d)''$. Since $\cla_0$ is also a 
type-I, strong mixing and Kolmogorov's property are equivalent. In particular the adjoint Markov semigroup 
$(\tilde{\tau}_t)$ on $\cla_0$ satisfies Kolmogorov's property, hence by Theorem 3.4,  
$||\tau_t -\phi_0|| \raro 0$ as $t \raro 0$. Hence by Theorem 3.3 $\bigcap_{n \ge 1} 
\Lambda_n(\pi_{\omega}(\otimes_{\IN}M_d)'')=\IC$. That is equivalent to 
$\{ S_IS^*_J: |I|=|J| < \infty \}''= \pi_{\omega}(\otimes_{\IN}M_d)''$. Since $S_i \in \pi_{\omega}(\otimes_{\IN}M_d)''$
by our construction we conclude that $\{ S_IS^*_J: |I|=|J| < \infty \}'' = \{ S_i,S^*_i: 1 \le i \le d \}''$. Thus 
for a type-I factor strongly mixing $(\cla_0,\eta_n,\psi_0)$ system, we find that $\cla_0=\clm$ 
and $l_k=V_k$.  \qed

\vsp
The above corollary enable us to construct a pure state on the UHF$_d$ algebra $\otimes_{\IZ}M_d$ so that its restriction 
on $\otimes_{\IN}M_d)$ is a type-I factor state (see [BJP], [BJKW], [Ma2]). For an explicit example of a pure state 
on $\otimes_{\IZ}M_d$ which give rise to a type-III factor on $\otimes_{\IN}M_d)$, we refer to [Ma1]. Is it possible 
to construct a pure state so that its restriction to one sided chain will be type-II? It is not hard to realize 
that it is impossible if we demand hyperfinite type-II$_1$ factor state. A proof and more results follow in the next 
section.

\newsection{ Jones index of a quantum dynamical semigroup on II$_1$ factor: }

In this section we continue our investigation in the general framework of section 4 and eventually study the case when 
$\cla_0$ is type-II$_1$.    

\vsp
\NI {\bf PROPOSITION 7.1: } Let $(\cla_0,\tau_t,\phi_0)$ be a dynamical system as in Theorem 4.1. If $\cla_{[0}$ is a 
type-II$_1$ factor which admits a unique normalize faithful normal tracial state then the following hold:

\NI (a) $F_{t]}=I$ for all $t \in \IR$;

\NI (b) $\tau=(\tau_t)$ is a semigroup of $*-$endomorphisms.

\NI (c) $\cla_{[0}=j_0(\cla_0)$.

\vsp
\NI {\bf PROOF:} Let $tr_0$ be the unique normalize faithful normal trace on $\cla_{[0}$. For any fix $t \ge 0$ we set a normal 
state $\phi_t$ on $\cla_{[0}$ by $\phi_t(x)=tr_0(\alpha_t(x))$. It is simple to check that it is also a faithful normal trace. 
Since $\alpha_t(I)=I$, by uniqueness $\phi_t=tr_0$. In particular $tr_0(F_{0]})=tr_0(\alpha_t(F_{0]})=tr_0(F_{t]})$, by faithful
property $F_{t]}=F_{0]}$ for all $t \ge 0$. Since $F_{t]} \uparrow 1$ as $t \raro \infty$ we have $F_{0]}=I$. 
Hence $F_{t]}=\alpha_t(F_{0]})=I$ for all $t \in \IR$. This proves (a). For (b) and (c) we recall that 
$F_{0]}j_t(x)F_{0]}=j_0(\tau_t(x))$ for all $t \ge 0$ and $j_t:\cla_0 \raro \cla_{[t}$ is an injective $*-$ homomorphism. 
Since $F_{t]}=F_{0]}=I$ we have $j_t(x)=F_{0]}j_t(x)F_{0]}= j_0(\tau_t(x))$. Hence $\cla_{[0}=j_0(\cla_0)$
and $j_0(\tau_t(x)\tau_t(y))=j_0(\tau_t(xy))$ for all $x,y \in \cla_0$. Now by injective property of $j_0$, 
we verify (b). This completes the proof. \qed

\vsp
\NI {\bf THEOREM 7.2:} Let $\omega'$ be a pure state on UHF$_d$ algebra $\otimes_{\IZ}M_d$ then the restriction $\omega$ of 
$\omega'$ to the UHF$_d$ algebra $\otimes_{\IN}M_d$ is not a hyperfinite type-II$_1$ factor. 

\vsp
\NI {\bf PROOF: } We first recall quantum dynamical semigroup $(\cla_0,\eta_t,\psi_0)$ described as in Theorem 6.4 and 
consider associated stationary minimal weak Markov processes $(j_t)$. The von-Neumann algebra $\cla_{[0}$ is isomorphic to 
a sub-algebra of $\pi(\otimes_{\IN}M_d)''$. In case $\pi(\otimes_{\IN}M_d)''$ is a type-II$_1$ factor, 
$\cla_0$ is also a type-II$_1$ factor. By Tomiyama's property [BR] we also note that $\cla_0$ is a hyperfinite factor
for hyperfinite $\pi(\otimes_{\IN}M_d)''$. In such a case by Theorem 4.1 $\cla_{[0}$ is also a type-II hyperfinite 
factor. Since $\cla_{[0}$ is isomorphic to a von-Neumann subalgebra of $\pi(\otimes_{\IN}M_d)''$ which is a 
type-II$_1$ factor, we conclude that $\cla_{[0}$ is also a hyperfinite type-II$_1$ factor. Hence by Proposition 7.1 (c), 
$(\eta_t)$ is a semigroup of $*-$endomorphisms and so $\psi_0(\eta_t(x)\eta_t(y))=\psi_0(xy)$ for all $t \ge 0$. Thus by 
Theorem 6.4 $\omega$ is pure if and only if $\psi_0(xy)=\psi_0(x)\psi_0(y)$ i.e. $x=\psi_0(x)I$ for all $x \in \cla_0$. This 
clearly contradicts that $\cla_0$ is a type-II$_1$ factor. \qed   

\vsp
We continue once more now our general case and fix a type-II$_1$ factor $\cla_0$ which admits a unique normalize faithful 
normal tracial state. Since $\cla_{[0}$ is a type-II factor whenever $\cla_0$ is so, we conclude that 
$\cla_{[0}$ is a type-II$_\infty$ factor whenever $\tau_t$ is not an endomorphism on a such a 
type-II$_1$ factor. The following proposition says much more.

\vsp
\NI {\bf THEOREM 7.3:}  Let $\cla_0$ be a type-II$_1$ factor with a unique normalize normal trace and $(\cla_0,\tau_t,\phi_0)$ 
be a dynamical system as in Theorem 4.1. Then the following hold:

\NI (a) $j_0(I)$ is a finite projection in $\cla_{[0}$ and there exists a type-II$_1$ factor $\clm_0$ isomorphic to $\cla_0$ so that 
$\cla_{[0}=\clm_0 \otimes \clb(\clf_{[0})$ where $\clf_{[0}$ is a complex separable Hilbert space. 

\NI (b) There exists a tower of type-II$_1$ factors $\clm_0 \subseteq \clm_s ...\subseteq \clm_t \subseteq ..,\;t \ge s \ge 0$ acting on 
a complex separable Hilbert space $\clf_0$ so that $\cla_{[-t}$ is isomorphic to $\clm_t \otimes \clb(\clf_{[t})$, where for each $t \ge 0$, 
$\clf_{[t}$ is a complex separable Hilbert space. 

\vsp
\NI {\bf PROOF:} By Theorem 4.1 $\cla_{[0}$ is a type-II factor. Thus $\cla_{[0}$ is either type-II$_1$ or type-II$_\infty$. In case it is 
type-II$_1$, Theorem 4.5 says that $\cla_{[-t}$ is $j_0(\cla_0)$, hence the statements (a) and (b) are true with $\clm_n=j_0(\cla_0)$ and
the Hilbert spaces $\clf_{[t}$ are $\IC$. Thus it is good enough if we prove (a) and (b) when $\cla_{[0}$ is indeed a type-II$_\infty$ 
factor. To that end we fix a normal faithful trace $tr$ on $\cla_{[0}$ and consider the map $x \raro tr(j_0(x)$ for 
$x \in \cla_0$. It is a normal faithful trace on $\cla_0$, hence it is a scaler multiple of the unique trace. Thus $j_0(I)$ is
a finite projection in $\cla_{[0}$. Now the general theory on von-Neumann algebra [Sak] guarantees the result once we recall that 
$j_0(\cla_0)=j_0(I)\cla_{[0}j_0(I)$. This proves the statement (a). For the second statement note also that $j_0(I)$ is a finite 
projection in $\cla_{[-t}$ for any $t \ge 0$, thus $j_0(I)\cla_{[-t}j_0(I)$ is a type-II$_1$ factor acting on $F_{0]}$. 
So once more we appeal to the general theory [Sak] for isomorphism with $\clm_t \otimes \clb(\clf_{[t})$. The 
inclusion relations follow from the inclusion relations $\cla_{[-s} \subseteq \cla_{[-t}$ where $t \ge s$ \qed

\vsp
We first recall Jones's index of a sub-factor originated to understand the structure 
of inclusions of von Neumann factors of type ${\rm II}_1$. Let $N$ be a sub-factor of a finite factor $M$. $M$ acts naturally as 
left multiplication on $L^2(M,tr)$, where $tr$ be the normalize normal trace. The projection $E_0=[N\omega] \in N'$, where $\omega$ 
is the unit trace vector i.e. $tr(x)=<\omega,x\omega>$ for $x \in M$, determines a conditional expectation $E(x)=E_0xE_0$ on $N$. If 
the commutant $N'$ is not a finite factor, we define the index $[M : N]$ to be infinite. In case $N'$ is also a finite factor, 
acting on $L^2(M,tr)$, then the index $[M : N]$ of sub-factors is defined as $tr(E_0)^{-1}$, which is the Murray-von Neumann 
coupling constant [MuN] of $N$ in the standard 
representation $L^2(M,tr)$. Clearly index is an invariance for the sub-factors. Jones proved $[M : N] \in \{4\cos^2(\pi/ n): n=3,4,
\cdots\}\cup[4,\infty]$ with all values being realized for some inclusion $N \subseteq M$.

\vsp
Let $\cla_0$ be a type-II$_1$ factor and $\tau$ be a normal completely positive unital normal map with a faithful normal invariant state 
$\phi_0$. We consider the dynamics $(\cla_0,\tau_n,\; n \ge 0,\phi_0)$ where $\tau_n=\tau \circ \tau  ...\circ \tau $ ($n$ fold ) and 
$\tau_0=I,$ the identity map and the associated tower of II$_1$ factors $\clm_0 \subseteq \clm_1 \subseteq ...\clm_n \subseteq $ acting 
on the complex separable Hilbert space $\clf_{0}$ described as in Theorem 4.6. Thus the infinite sequence of Jones index 
$\{ [\clm_{k+1} : \clm_{k}],\; k=0,1,..\} $ associated with the canonical tower $\clm_k: k \ge 0 $ of II$_1$ factors is an 
invariance for the dynamics $(\cla_0,\tau,\phi_0)$. One aim to investigate how this tower $(\clm_k: \ge 0)$ is related with Jones's
tower of type-II$_1$ factors.   

\vsp
To that end we review now Jones's construction [Jo, OhP]. Let $\phi_0$ be the unique normalize normal trace. 
The algebra $\cla_0$ acts on $L^2(\cla_0,\phi_0)$ by left multiplication $\pi_0(y)x=yx$ for $x \in L^2(\cla_0,\phi_0)$. Let $\omega$ be the cyclic and separating trace vector in $L^2(\cla_0,\phi_0)$. The projection $E_0=[\clb_0\omega]$ induces a trace preserving conditional 
expectation $\tau: a \raro E_0aE_0$ of $\cla_0$ onto $\clb_0$. Thus $E_0\pi_0(y)E_0=E_0\pi_0(E(y))E_0$ for all $y \in \cla_0$. Let 
$\cla_1$ be the von-Neumann algebra $\{ \pi_0(\cla_0),E_0 \}''$. $\cla_1$ is also a type-II$_1$ factor and 
$\cla_0 \subseteq \cla_1$, where we have identified $\pi_0(\cla_0)$ with $\cla_0$. Jones proved that $[\cla_1 : \cla_0]=[\cla_0 : \clb_0]$. Now by repeating this canonical method we get an increasing tower of type-II$_1$ factors $\cla_1 \subseteq \cla_2 ...$ so that 
$[\cla_{k+1} : \cla_k]=[\cla_0 : \clb_0]$ for all $k \ge 0$. Thus the natural question: How Jones tower 
$\cla_0 \subseteq \cla_1 \subseteq ... \subseteq \cla_k ...$ is related 
with the tower $\clm_0 \subseteq \clm_1 ... \clm_k \subseteq \clm_{k+1}$ associated with the dynamics $(\cla_0,\tau_n,\phi_0)$?     

\vsp
To that end recall the von-Neumann sub-factors $\clm_0 \subseteq \clm_1$ and the induced representation of $\clm_1$ 
on Hilbert subspace $H_{[-1,0]}$ generated by $\{ j_0(x_0)j_{-1}(x_{-1})\Omega: x_0,x_{-1} \in \cla_0 \}$. $\Omega$ is 
the trace vector for $\clm_0$ i.e. $\phi_0(x)=<\Omega,j_0(x)\Omega>$. It is the trace vector for $\clm_1$ if and only if 
$\clm_1=\clm_0$, ( for trace vector we check that $\phi_0(\tau(x)y\tau(z)=\phi(j_0(x)j_{-1}(y)j_0(z) = \phi_0(\tau(zx)y)$) 
for any $x,y,z \in \cla_0$ ). Nevertheless there exists a unique normalize trace on $\clm_1$, being a type-II$_1$ factor. 

\vsp
\NI {\bf THEOREM 7.4: } $[\clm_1 : \clm_0] > [\cla_0 : \clb_0]$.  

\vsp
\NI {\bf PROOF: } Let $\phi_1$ be the unique normalize normal trace on $\cla_1$ and $\clh_1 = L^2(\cla_1, \phi_1)$. We consider 
the left action $\pi_1(x):y \raro xy$ of $\cla_1$ on $\clh_1$. Thus $\pi_0(\cla_0)$ is also acting on $\clh_1$. Since $E_0\pi_0(x)E_0
=E_0\pi_0(\tau(x))E_0=E_0\pi_0(\tau(x))$, for any element $X \in \cla_1$, $E_0X=E_0\pi_0(x)$ for some $x \in \cla_0$. Thus $\pi_1(E_0)$ is the projection on the subspace $\{ E_0\pi_0(x): x \in \cla_0 \}$.  

\vsp
For any $y \in \cla_0$ we set

\NI (a) $k_{-1}(y)$ on the subspace $\pi_1(E_0)$ by $j_{-1}(y)E_0\pi_0(x)=E_0\pi_0(yx)$ for $x \in \cla_0$ and extend it to 
$\clh_1$ trivially.   

\NI (b) $k_0(y)x =\pi_0(y)x$ for $x \in \cla_1$

\vsp
For $y,z \in \cla_0$ we verify that 
$$<E_0\pi_0(y), k_{-1}(1)k_0(x)k_{-1}(1)E_0\pi_0(z)>_1 = <E_0\pi_0(y),E_0\pi_0(x)E_0\pi_0(z)>_1 $$ 
$$= <E_0\pi_0(y),E_0\pi_0(\tau(x))E_0\pi_0(z)>_1$$
$$= <E_0\pi_0(y),E_0\pi_0(\tau(x))\pi_0(z)>_1$$
Thus $k_{-1}(1)k_0(x)k_{-1}(1)=k_{-1}(\tau(x))$ for all $x \in \cla_0$. Note that $k_{-1}(1)=\pi_1(E_0)$ and the identity 
operator in $\clh_1$ is a cyclic vector for the weak Markov process and thus by uniqueness of minimal weak Markov processes 
associated with $(\cla_0,\tau,\phi_0)$, $\{ k_{-1}(\cla_0), k_0(\cla_0) \}''$ is isomorphic to $\clm_1$. 
Since $k_{-1}(1)=\pi_1(E_0)$, $\cla_1 \subseteq \clm_1$. In fact strict inclusion hold unless $\clb_0=\cla_0$. Thus 
$[\clm_1 : \cla_1] > 1$. Since $[\clm_1 : \clm_0]= [\clm_1 : \cla_1] [\cla_1 : \clm_0]$ and $[\cla_1 : \cla_0] = 
[\cla_0 : \clb_0]$, we conclude the result. \qed

\vsp
Thus for any finite sub-factor $\clb_0$ of a type-II$_1$ factor we could associated via a canonical method a sequence
of Jones numbers $\{ [\clm_k : \clm_{k-1}],\; k \ge 0 \}$ as invariance for the inclusion $\clb_0 \subseteq \cla_0$, 
where we set $\clm_{-1}=\clb_0$.

\bigskip
{\centerline {\bf REFERENCES}}

\begin{itemize} 

\bigskip
\item{[AM]} Accardi, L., Mohari, A.: Time reflected Markov processes. Infin. Dimens. Anal. Quantum Probab. Relat. Top., vol-2, no-3, 397-425 (1999).

\item{[AW]} Araki, H., Woods, E.J.: Complete boolean algebras of type I factors, Publ. Res. Inst. Math. Sci. Series A, vol-II (1966) 157-242. 

\item {[Ar]} Arveson, W.: Pure $E_0$-semigroups and absorbing states, Comm. Math. Phys. 187 , no.1, 19-43, (1997)

\item {[Bh]} Bhat, B.V.R.: An index theory for quantum dynamical semigroups, Trans. Amer. Maths. Soc. vol-348, no-2 561-583 (1996).   

\item {[BP]} Bhat, B.V.R., Parthasarathy, K.R.: Kolmogorov's existence theorem for Markov processes on $C^*$-algebras, Proc.
Indian Acad. Sci. 104,1994, p-253-262.

\item {[BR]} Bratelli, Ola., Robinson, D.W. : Operator algebras and quantum statistical mechanics, I,II, Springer 1981.

\item {[BJ]} Bratteli, Ola; Jorgensen, Palle E. T. Endomorphism of $\clb(\clh)$, II,
Finitely correlated states on $\clo_N$, J. Functional Analysis 145, 323-373 (1997).

\item {[BJP]} Bratelli, Ola., Jorgensen, Palle E.T. and Price, G.L.: Endomorphism of $\clb(\clh)$, Quantization, nonlinear partial differential 
equations, Operator algebras, ( Cambridge, MA, 1994), 93-138, Proc. Sympos. Pure Math 59, Amer. Math. Soc. Providence, RT 1996.

\item {[BJKW]} Bratelli, O., Jorgensen, Palle E.T., Kishimoto, Akitaka and
Werner Reinhard F.: Pure states on $\clo_d$, J.Operator Theory 43 (2000),
no-1, 97-143.

\item{[Da]} Davies, E.B.: Quantum Theory of open systems, Academic press, 1976.

\item{[El]} Elliot. G. A.: On approximately finite dimensional von-Neumann algebras I and II, Math. Scand. 39 (1976), 91-101;
Canad. Math. Bull. 21 (1978), no. 4, 415--418. 

\item{[FNW1]} Fannes, M., Nachtergaele,D., Werner,R.: Finitely Correlated States on Quantum Spin Chains, Commun. Math. Phys. 144,
443-490 (1992).

\item{[FNW2]} Fannes, M., Nachtergaele,D., Werner,R.: Finitely Correlated pure states, J. Funct. Anal. 120, 511-534
(1994).

\item{[Fr]} Frigerio, A.: Stationary states of quantum dynamical
semigroups. Commun. Math. Phys. 63, 269-276 (1978).

\item{[Jo]} Jones, V. F. R.: Index for subfactors. Invent. Math. 72 (1983), no. 1, 1--25.

\item{[Li]} Lindblad, G. : On the generators of quantum
dynamical semigroups, Commun.  Math. Phys. 48, 119-130 (1976).

\item{[Ma1]} Matsui, T.: A characterization of pure finitely correlated states.
Infin. Dimens. Anal. Quantum Probab. Relat. Top. 1 (1998), no. 4, 647--661.

\item{[Ma2]} Matsui, T.: The split property and the symmetry breaking of the quantum spin chain, Comm.
Maths. Phys vol-218, 293-416 (2001).

\item{[Mo1]} Mohari, A.: Markov shift in non-commutative probability, Jour. Func. Anal. 199 (2003) 189-209.  

\item{[Mo2]} Mohari, A.: Endomorphisms on $\clb(\clh)$, preprint.

\item{[Mo3]} Mohari, A.: Quantum detailed balance and split property in quantum spin chain, Submitted to Communication in
Mathematical Physics, 2004.
 
\item{[MuN]} Murray, F. J.; von Neumann, J., On rings of operators. (English)[J] Ann. Math., Princeton, (2)37, 116-229.

\item{[OP]} Ohya, M., Petz, D.: Quantum entropy and its use, Text and monograph in physics, Springer-Verlag 1995. 

\item{[Po]} Powers, Robert T.: An index theory for semigroups of $*$-endomorphisms of
$\clb(\clh)$ and type II$_1$ factors.  Canad. J. Math. 40 (1988), no. 1, 86--114.

\item{[Sak]} Sakai, S.: C$^*$-algebras and W$^*$-algebras, Springer 1971.  

\item{[Sa]} Sauvageot, Jean-Luc: Markov quantum semigroups admit covariant Markov $C^*$-dilations. Comm. Math. Phys. 
106 (1986), no. 1, 91­103.

\item{[So]} Stormer, Erling : On projection maps of von Neumann algebras. Math. Scand. 30 (1972), 46--50.

\item{[Vi]} Vincent-Smith, G. F.: Dilation of a dissipative quantum dynamical system to a quantum Markov process. Proc. 
London Math. Soc. (3) 49 (1984), no. 1, 58­72. 

\end{itemize}

\end{document}